\documentclass{cocv}
\usepackage{amsmath,amsthm,a4wide}

  {\par\medbreak\refstepcounter{theorem}%
    \noindent\textbf{Remark~\thetheorem. }}%
  {\par\medskip}

\newcommand{\vz}[1]{\ensuremath{\mathbb{#1}}}
\newcommand{\R}{{\vz R}}
\newcommand{\N}{{\vz N}}

\def\ubc{u_{\mathrm{bc}}}
\let\upGamma\Gamma
\def\Gamma{\mathit\upGamma}

\def\pref#1{(\ref{#1})}

\long\def\drop#1{}
\def\longrightharpoonup{\relbar\joinrel\rightharpoonup}
\let\weakto\longrightharpoonup

\let\downto\downarrow
\let\e\varepsilon
\let\epsilon\varepsilon

\def\T{\vz T}

\let\ds\displaystyle
\def\div{\mathop{\mathrm{div}}}
\def\XXint#1#2#3{{\setbox0=\hbox{$#1{#2#3}{\int}$}
     \vcenter{\hbox{$#2#3$}}\kern-.5\wd0}}

\usepackage{latexsym}
\usepackage{amssymb}
\usepackage{amsfonts}
\usepackage{bm}
\usepackage{epsfig}
\usepackage{bbm}
\usepackage{amsmath}
\usepackage{graphicx}
\usepackage{psfrag}
\usepackage{color}
\usepackage{multicol}
\usepackage{subfig}
\usepackage{appendix}

\begin{document}
%%-----------------------------
%%      the top matter
%%-----------------------------
\title{The $H^{-1}$-norm of tubular neighbourhoods of curves}

\author{Yves van Gennip}\address{Department of Mathematics Simon Fraser University, 8888 University Drive, Burnaby, British Columbia  V5A 1S6, Canada}
\author{Mark A. Peletier}\address{Dept. of Mathematics and Computer Science, Technische Universiteit Eindhoven, PO Box 513, 5600 MB  Eindhoven, The Netherlands}
\date{}
\begin{abstract}
We study the $H^{{-1}}$-norm of the function $1$ on tubular neighbourhoods of curves in $\R^{2}$. We take the limit of small thickness $\e$, and we prove two different asymptotic results. The first is an asymptotic development for a fixed curve in the limit $\e\to0$, containing contributions from the length of the curve (at order $\e^{3}$), the ends ($\e^{4}$), and the curvature ($\e^{5}$).

The second result is a $\upGamma$-convergence result, in which the central curve may vary along the sequence $\e\to0$. We prove that a rescaled version of the $H^{-1}$-norm, which focuses on the $\e^{5}$ curvature term, $\upGamma$-converges to the $L^2$-norm of curvature. In addition, sequences along which the rescaled norm is bounded are compact in the $W^{1,2}$-topology. 

Our main tools are the maximum principle for elliptic equations and the use of appropriate trial functions in the variational characterisation of the $H^{-1}$-norm. For the $\upGamma$-convergence result we use the theory of systems of curves without transverse crossings to handle potential intersections in the limit. 
\end{abstract}
\begin{resume} 
Nous \'etudions la norme $H^{-1}$ de la fonction $1$ sur des domaines minces dans~$\R^2$. Nous consid\'erons des suites  de voisinages  tubulaires de courbes planes. Nous d\'emon\-trons deux caract\'erisations asymptotiques de cette norme dans la limite de petite largeur $\e$.

Le premier r\'esultat est un d\'eveloppement asymptotique pour les $\e$-voisinages tubulaires d'une courbe fixe. Dans ce d\'eveloppement apparaissent des termes provenant de la longueur de la courbe (\`a l'ordre $\e^3$), des extremit\'es ($\e^4$) et de la courbure ($\e^5$).

Le deuxi\`eme r\'esultat concerne des suites d'$\e$-voisinages de courbes, dans le cas o\`u les courbes peuvent varier le long de la suite. Nous d\'emontrons que la norme $H^{-1}$ $\upGamma$-converge vers la norme $L^2$ de la courbure.  Cette $\upGamma$-convergence a lieu par rapport \`a la topologie $W^{1,2}$, et une suite dont la norme renormalis\'ee est born\'ee est compacte dans cette topologie.

Les preuves font appel au principe du maximum pour les \'equations elliptiques et \`a une caract\'erisation variationnelle de la norme $H^{-1}$. Pour la $\upGamma$-convergence, la th\'eorie de syst\`emes de courbes sans intersections transverses permet de traiter les intersections dans la limite. 

\end{resume}
\maketitle
%%-----------------------------
%%      your text
%%-----------------------------

\section{Introduction}
In this paper we study the set function $F:2^{\R^2}\to\R$, 
\[
F(\Omega) := \|1\|_{H^{-1}(\Omega)}^2 := \sup \left\{ \int_\Omega \left[ 2 u - |\nabla u|^2\right]\,dx : u \in C_c^{\infty}(\Omega)\right\}.
\]
More specifically, we are interested in the value of $F$ on \emph{$\e$-tubular neighbourhoods} $T_\e\gamma$ of a curve~$\gamma$, i.e. on the set
of points strictly within a distance $\e$ of $\gamma$.

The aim of this paper is to explore the connection between the geometry of a curve~$\gamma$ and the values of $F$ on the $\e$-tubular neighbourhood $T_\e\gamma$. Our first main result is the following asymptotic development. If $\gamma$ is a smooth open curve, then
\begin{equation}
\label{devel:main}
\|1\|_{H^{-1}(T_\e\gamma)}^2 = \frac23 \e^3 \ell(\gamma) + 2\alpha \e^4 + 
  \frac 2{45}\e^5\int_{\gamma} \kappa^2 + O(\e^6)\qquad\text{as }\e\to0.
\end{equation}
Here $\ell(\gamma)$ is the length of $\gamma$, $\alpha>0$ is a constant independent of $\gamma$, and $\kappa$ is the curvature of~$\gamma$. The `2' that multiplies $\alpha$ in the formula above is actually the number of end points of $\gamma$; for a closed curve the formula holds without this term. Under some technical restrictions~\pref{devel:main} is proved in Theorem~\ref{th:devel}.

The expansion~\pref{devel:main} suggests that for \emph{closed} curves the rescaled functional 
\[
G_\e(\gamma):= \e^{-5} \left(\|1\|_{H^{-1}(T_\e\gamma)}^2 - \frac23 \e^3\ell(\gamma)\right)
\]
resembles the \emph{elastica functional}
\[
G_0(\gamma) := \frac2{45}\int_\gamma \kappa^2.
\]
With our second main result we convert this suggestion into a $\upGamma$-convergence result, and supplement it with a statement of compactness. Before we describe this second result in more detail, we first explain the origin and relevance of this problem.

\subsection{Motivation}\label{sec:motivation}

The $H^{-1}$-norm of a set or a function appears naturally in a number of applications, such as electrostatic interaction or gravitational collapse. The case of tubular neighbourhoods and the relationship with geometry are more specific. We mention two different origins.

The discussion of the connection between the geometry of a domain and the eigenvalues of the Laplacian goes back at least to H. A. Lorentz' Wolfskehl lecture in 1910, and has been popularized by Kac's and Bers' famous question `can one hear the shape of a drum?'~\cite{Kac66}. The first eigenvalue of the Laplacian with Dirichlet boundary conditions is actually strongly connected to the $H^{-1}$-norm. 
This relation can be best appreciated when writing the definition of the first eigenvalue under Dirichlet boundary conditions as
\begin{equation}
\label{def:lambda0}
\lambda_0(\Omega) = \inf \;\left\{\frac{\ds\int_{\Omega} |\nabla u|^2 }{\ds\int_{\Omega} u^2 }: 
u\in C_c^\infty(\Omega)\right\},
\end{equation}
and the $H^{-1}$-norm as
\begin{equation}
\label{def:HMO2}
\|1\|_{H^{-1}(\Omega)}^2 = \sup\; \left\{\frac{\ds\Bigl(\int_\Omega u\Bigr)^2}{\ds\int_\Omega |\nabla u|^2}:
u\in C_c^\infty(\Omega)\right\}.
\end{equation}
Sidorova and Wittich~\cite{SidorovaWittich08} investigate the $\e$- and $\gamma$-dependence of $\lambda_0(T_\e\gamma)$. As in the case of the $H^{-1}$-norm, the highest-order behaviour of $\lambda_0(T_\e\gamma)$ is dominated by the short length scale~$\e$ alone; the correction, at an order $\e^2$ higher, depends on the square curvature. The signs of the two correction terms are different, however: while the curvature correction in $\|1\|^2_{H^{-1}}$ (the third term on the right-hand side of~\pref{devel:main}) comes with a positive sign, this correction carries a negative sign in the development of $\lambda_0$. 

This sign difference can also be understood from the difference between~\pref{def:lambda0} and~\pref{def:HMO2}. 
Assume that for a closed curve the supremum in~\pref{def:HMO2} is attained by $\hat u$. The development in~\pref{devel:main} states that for small $\e$,  $\bigl(\int_{T_\e\gamma} \hat u\bigr)^2 / \int_{T_\e\gamma} |\nabla \hat u|^2 \approx \e^3 C_1 (1+ \e^2 C_2)$, for two positive constants $C_1$ and $C_2$ that depend only on the curve. Inverting the ratio, we find that $\int_{T_\e\gamma} |\nabla \hat u|^2 / \bigl(\int_{T_\e\gamma} \hat u\bigr)^2 \approx \e^{-3} C_1^{-1}  (1-\e^2 C_2)$. If we disregard the distinction between $\int u^2$ and $(\int u)^2$, then this argument explains why the curvature correction enters with different signs.

\medskip
The question that originally sparked this investigation was that of partial localisation.
Partial localisation is a property of certain pattern-forming systems. The term `localisation' refers to structures---e.g. local or global energy minimisers---with limited spatial extent. `Partial localisation' refers to a specific subclass of structures, which are localised in some directions and extended in others. Most systems tend to either localise in all directions, such as in graviational collapse, or to delocalise and spread in all directions, as in diffusion. {Stable} partial localisation is therefore a relatively rare phenomenon, and only a few systems are known to exhibit it~\cite{PeletierRoeger08,DAprile00,DoelmanVanderPloeg02, vanGennipPeletier08, vanGennipPeletier09a}

In two dimensions, partially localised structures appear as fattened curves, or when their boundaries are sharp, as tubular neigbourhoods. Previous work of the authors suggests that various energy functionals all involving the $H^{-1}$-norm might exhibit such partial localisation, and some existence and stability results are already available~\cite{vanGennipPeletier08,vanGennipPeletier09a}. On the other hand the partially localising property of these functionals without restrictions on geometry is currently only conjectured, not proven. The work of this paper can be read as an intermediate step, in which the geometry is partially fixed, by imposing the structure of a tubular neighbourhood, and partially free, by allowing the curve $\gamma$ to vary.

The freedom of variation in $\gamma$ gives rise to questions that go further than a simple asymptotic development in $\e$ for fixed $\gamma$. A common choice in this situation is the concept of $\upGamma$-convergence; this concept of convergence of functionals implies convergence of minimisers to minimisers, and is well suited for asymptotic analysis of variational problems. For this reason our second main result is on the $\upGamma$-convergence of the functional $G_\e$.

\medskip

Before we state this result in full, we first comment on curvature and regularity, and we then introduce the concept of \emph{systems of curves}.

\subsection{Curvature and regularity}
\label{sec:intro-curvature}

In this paper we only consider the case in which the tubular neighbourhoods are regular, in the following sense, at least for sufficiently small $\e$: for each $x\in T_\e\gamma$
there exists a unique point $\tilde x\in (\gamma)$ of minimal distance to $x$, where $(\gamma)\subset\R^2$ is the trace or image of the curve $\gamma$.
An equivalent formulation of this property is given in terms of an upper bound on the \emph{global radius of curvature} of $\gamma$:

\begin{dfntn}[\cite{GonzalezMaddocks99}]\label{def:globalcurvature}
If $x, y, z\in (\gamma)$ are pairwise disjoint and not collinear, let $r(x, y, z)$ be the radius of the unique circle in $\R^2$ through $x$, $y$, and $z$ (and let $r(x,y,z)=\infty$ otherwise). The \emph{global radius of curvature} of $\gamma$ is defined as
\[
\rho(\gamma) := \underset{x, y, z \in (\gamma)}{\inf} r(x, y, z).
\]
\end{dfntn}
Since the `local' curvature $\kappa$ is bounded by $1/\rho(\gamma)$, finiteness of the global curvature implies $W^{2,\infty}$-regularity of the curve. More specifically, regularity of the $\e$-tubular neighbourhood $T_\e\gamma$ is equivalent to the statement $\rho(\gamma)\geq \e$.

\subsection{Systems of curves}

\begin{figure}[h]
\hspace{0.025\textwidth}
\subfloat[][Curve with locally multiplicity 2, limit of a sequence of single smooth simple curves]
{
    \includegraphics[height=4cm]{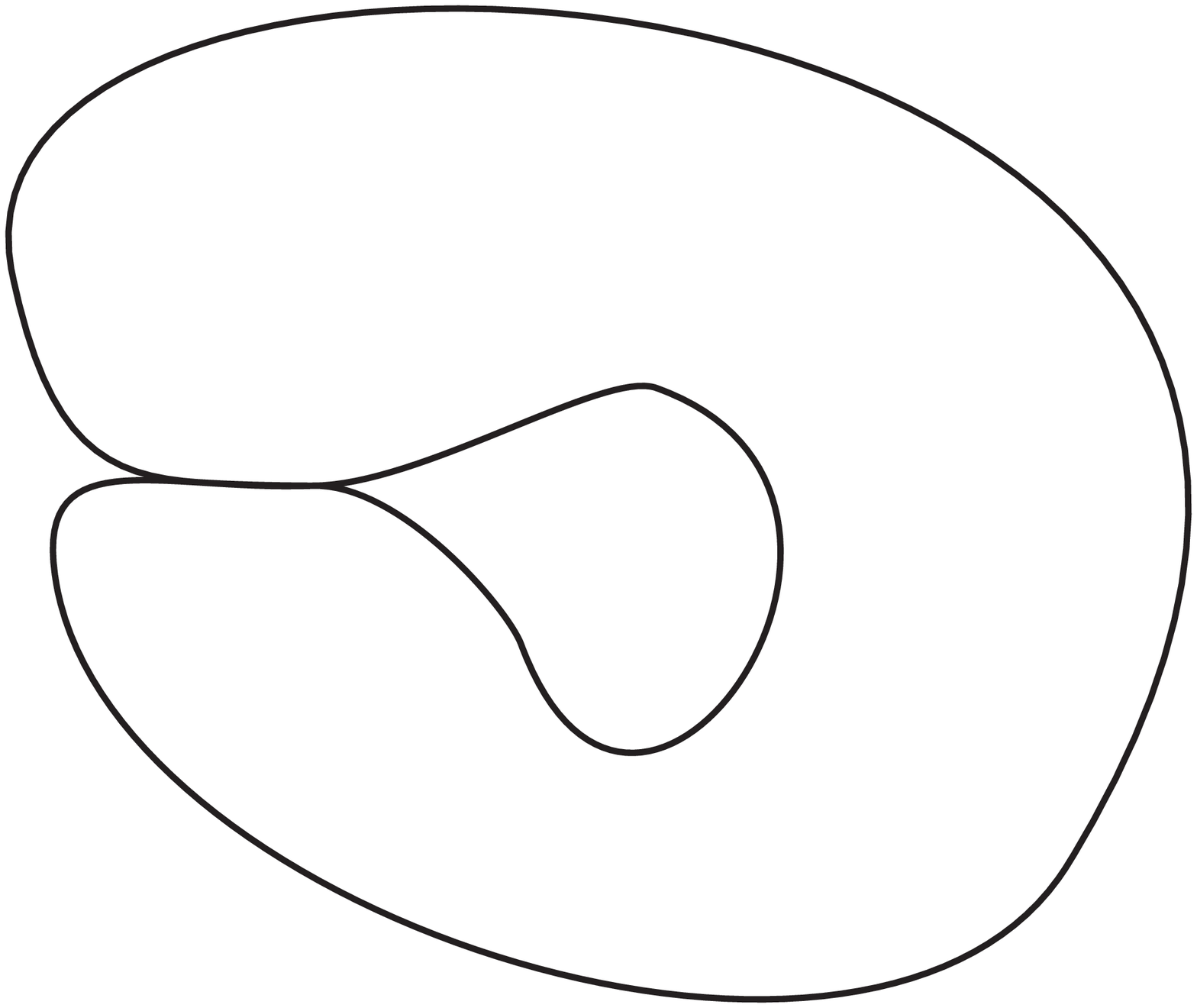}
    \label{fig:whysystems}	
}
\hspace{0.2\textwidth}
\subfloat[whysystems][This curve is not a limit of single simple curves, but can be obtained as the limit of a sequence of pairs of simple curves]
{
    \qquad\qquad\includegraphics[width=2cm]{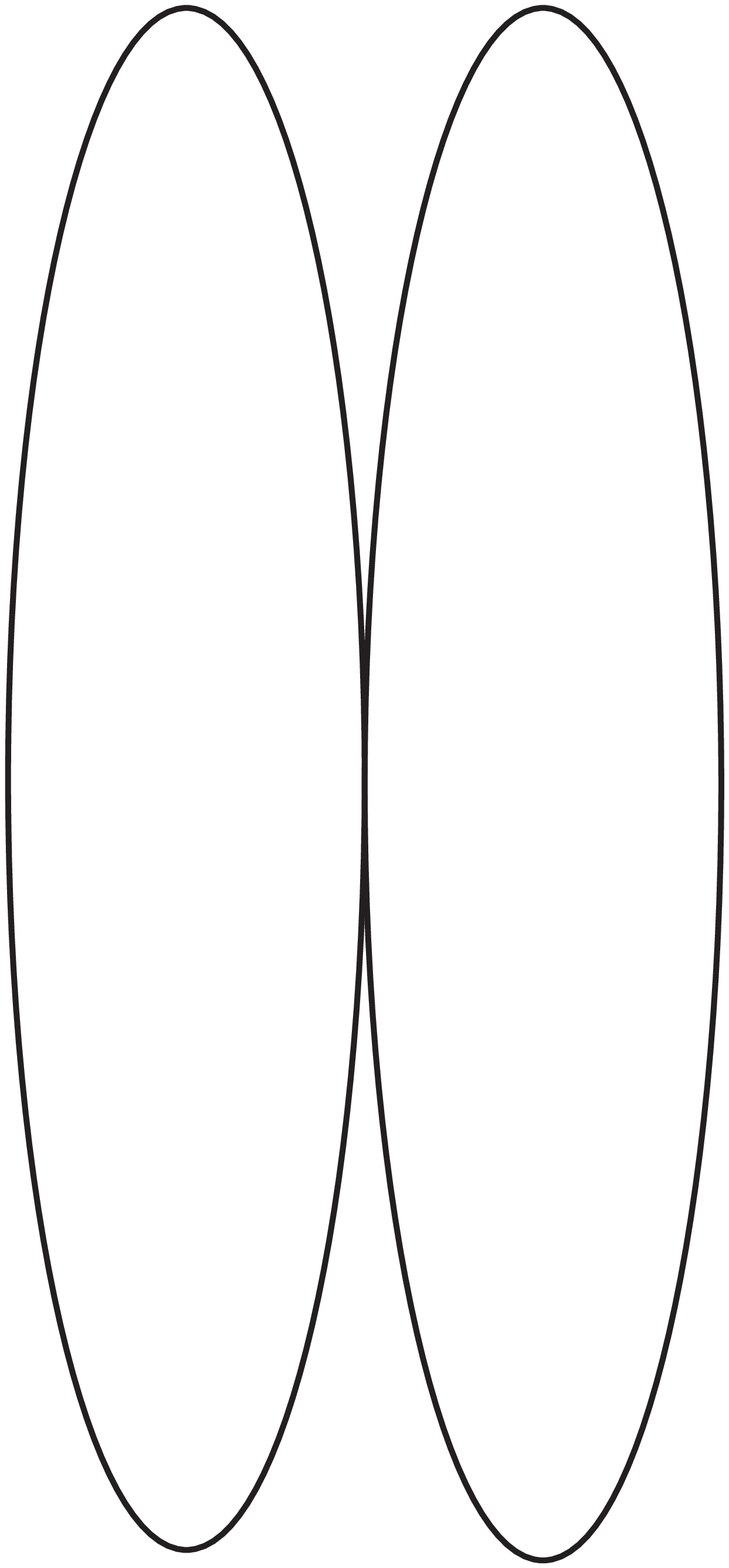} \qquad \qquad \qquad
    \label{fig:whysystems2}
}
\caption{Two curves with locally multiplicity 2.}
\label{fig:non-simple}
\end{figure}

Neither compactness nor $\upGamma$-convergence of $G_\e$ is expected to hold for simple, smooth closed curves, where `simple' means `non-self-intersecting'. One reason is that a perfectly reasonable sequence of simple smooth closed curves may converge to a non-simple curve, as shown in Figure~\ref{fig:whysystems}. 
Nothing in the energy $G_\e$ will prevent this; therefore we need to consider a generalisation of the concept of a simple closed curve.

The work of Bellettini and Mugnai~\cite{BellettiniMugnai04,BellettiniMugnai07}
 provides the appropriate concept. Leaving aside issues of regularity for the moment (the full definition is given in Section~\ref{sec:systemsclosedcurves}), a \emph{system of curves without transverse crossings} $\Gamma$ is a finite collection of curves,  $\Gamma=\{\gamma_i\}_{i=1}^m$, with the restriction that 
\[
\gamma_i(s) = \gamma_j(t) \quad \text{for some }i,j,s,t \qquad\Longrightarrow \qquad
\gamma_i'(s) \parallel \gamma_j'(t).
\]
In words: intersections are allowed, but only if they are tangent. 
Continuing the convention for curves, we write $(\Gamma)$ for the trace of $\Gamma$, i.e. $(\Gamma):= \bigcup_{i=1}^m (\gamma_i)\subset\R^2$. The \emph{multiplicity} $\theta$ of any point $x\in(\Gamma)$ is given by
\[
\theta(x) := \#\{(i,s): \gamma_i(s) = x\}.
\]
Figure~\ref{fig:whysystems} is covered by this definition, by letting $\Gamma$ consist of a single curve $\gamma$, and where $\theta$ equals $2$ on the intersection region and $1$ on the rest of the curve. 

Figure~\ref{fig:whysystems2} is an example of a system of curves without transverse crossings which can be represented by either one or two curves $\gamma$. This example motivates the introduction of an equivalence relationship on the collection of such systems. Two systems of curves $\Gamma^1$ and $\Gamma^2$ are called \emph{equivalent} if $(\Gamma^1) = (\Gamma^2)$ and $\theta^1 \equiv \theta^2$; this relationship gives rise to \emph{equivalence classes} of such systems of curves without transverse crossings.

\medskip

This leads to the definition of the sets $SC^{1,2}$ and $SC^{2,2}$, whose elements are \emph{equivalence classes of systems of curves}, for which each curve is of regularity $W^{1,2}$ or $W^{2,2}$. All admissible objects will actually be elements of $SC^{2,2}$; the main use of $SC^{1,2}$ is to provide the
right concept of convergence in which to formulate the compactness and $\upGamma$-convergence below. Where necessary, we write $[\Gamma]$ for the equivalence class (the element of $SC^{k,2}$) containing $\Gamma$; where possible, we simply write $\Gamma$ to alleviate notation.

\subsection{Compactness and $\upGamma$-convergence}

With this preparation we can state the second main result of this paper. The discussion above motivates changing the definition of the functionals $G_\e$ and $G_0$ defined earlier to incorporate conditions on global curvature and to allow for systems of curves. Note that in this section we only consider systems of \emph{closed} curves.

Define the functional $\mathcal{G}_{\epsilon}: SC^{1,2} \to \R\cup\{\infty\}$ by
\[
\mathcal{G}_{\epsilon}(\Gamma) := 
  \left\{\begin{array}{ll}
     \epsilon^{-5} \|1\|_{H^{-1}(T_\e\Gamma)}^2-\frac23\epsilon^{-2} \ell(\Gamma)
       & \text{ if }\Gamma\in SC^{2,2} \text{ and }\rho(\Gamma)\geq \e\\
     +\infty&\text{ otherwise},
  \end{array}\right.
\]
and let $\mathcal{G}_0: SC^{1,2} \to [0, \infty]$ be defined by
\[
\mathcal{G}_0(\Gamma) := 
 \left\{\begin{array}{ll}
  \displaystyle\frac2{45}\sum_{i=1}^m\ell(\gamma_i) \int_{\gamma_i} \kappa_i^2 
      & \text{ if }\Gamma\in\mathcal{S}_0,\\
  +\infty&\text{ otherwise,}
 \end{array}\right.
\]
where $\Gamma = \{\gamma_i\}_{i=1}^m$ and $\kappa_i$ is the curvature of $\gamma_i$, and where the admissible set $\mathcal S_0$ is given by
\begin{align*}
\mathcal{S}_0&:= \Bigl\{ \Gamma\in SC^{2,2}: \ell(\Gamma)<\infty \text{ and } \Gamma \text{ has no transverse crossings} \Bigr\}.
\end{align*}

The values of $\mathcal{G}_{\epsilon}(\Gamma)$ and $\mathcal{G}_0(\Gamma)$ are
independent of the choice of representative (see Remark~\ref{rem:indepofrep}), so that $\mathcal{G}_{\epsilon}$ and $\mathcal{G}_0$ are well-defined on equivalence classes. 

\medskip
We have compactness of energy-bounded sequences, provided they have bounded length and remain inside a fixed bounded set:
\begin{thrm}
\label{th:compactness}
Let $\e_n\downto0$, and let $\{\Gamma^n\}_{n\geq 1}\subset SC^{2,2}$ be a sequence such that 
\begin{itemize}
\item There exists $R>0$ such that $(\Gamma^n)\subset B(0,R)$ for all $n$;
\item $\sup_n \ell(\Gamma^n) < \infty$;
\item $\sup_n \mathcal{G}_{\e_n}(\Gamma^n)< \infty$.
\end{itemize} 
Then $\Gamma^n$ converges along a subsequence to a limit $\Gamma\in \mathcal S_0$ in the convergence of $SC^{1,2}$.
\end{thrm}
The concept of convergence in $SC^{1,2}$ is defined in Section~\ref{sec:systemsclosedcurves}. In addition to this compactness result, the functional $\mathcal G_0$ is the $\upGamma$-limit of $\mathcal G_\e$:
\begin{thrm}
\label{thm:gammaconvergence}
Let $\e_n\downto 0$. 
\begin{enumerate}
\item 
\label{thm:gammaconvergence:1}
If $\Gamma_n\in SC^{1,2}$ converges to $\Gamma\in SC^{1,2}$
in the convergence of $SC^{1,2}$, then $\mathcal{G}_0(\Gamma) \leq \underset{n\to\infty}{\liminf}\, \mathcal{G}_{\epsilon_n}(\Gamma^n)$.
\item \label{thm:gammaconvergence:2}
If $\Gamma \in SC^{1,2}$, then there is a sequence $\{\Gamma^n\}_{n\geq1}\subset SC^{1,2}$ converging to $\Gamma$
in the convergence of convergence of $SC^{1,2}$ for which $\mathcal{G}_0(\Gamma) \geq \underset{n\to\infty}{\limsup}\, \mathcal{G}_{\epsilon_n}(\Gamma^n)$.
\end{enumerate}
\end{thrm}

\subsection{Discussion}

\subsubsection{Hutchinson varifolds}

There is a close relationship between the systems of curves of Bellettini \& Mugnai and a class of varifolds. 
To a system of curves $\Gamma:=\{\gamma_i\}_{i=1}^m$ we can associate a measure $\mu_{\Gamma}$ via
\[
\int_{\R^2} \varphi\,d\mu_{\Gamma} = \sum_{i=1}^m \int_\gamma \varphi(\gamma_i(s)) |\gamma_i'(s)|\,ds,
\]
for all $\varphi \in C_c(\R^2)$. By \cite[Remark 3.9, Proposition 4.7, Corollary 4.10]{BellettiniMugnai07}  $\Gamma$ is a $W^{2,2}$-system of curves without transverse crossings if and only if $\mu_{\Gamma}$ is a Hutchinson varifold (also called curvature varifold) with weak mean curvature
$H\in L^2(\mu_{\Gamma})$,
such that a unique tangent line exists in every $x\in(\Gamma)$. Two systems of curves are mapped to the same varifold if and only if they are equivalent, a property which underlines that the appropriate object of study is the equivalence class rather than the system itself.

The compactness result for integral varifolds (\cite[Theorem 6.4]{Allard72}, \cite[Theorem 3.1]{Hutchinson86}) can be extended to a result for Hutchinson varifolds under stricter conditions which imply a uniform control on the second fundamental form along the sequence (\cite[Theorem 5.3.2]{Hutchinson86}). In our case we do not have such a control on the curvature, since the bound on the global radius of curvature, $\rho(\Gamma)\geq \epsilon$ vanishes in the limit  $\epsilon\to 0$. Therefore the compactness result of Theorem~\ref{th:compactness} covers a situation not treated by Hutchinson's result.

\subsubsection{Extensions}

The current work opens the way for many extensions that can serve as the subject of future inquiries. 
One such is the proof of a $\upGamma$-convergence result that also takes open curves into account. Expansion (\ref{devel:main}) suggests two possible functionals for study: 
\[
\e^{-4} \left(\|1\|_{H^{-1}(T_\e\Gamma)}^2 - \frac23 \e^3\ell(\Gamma)\right),
\]
which is expected to approximate
\[
2 n(\Gamma) \alpha,
\]
where $n(\Gamma)$ is the number of open curves in $\Gamma$; the second functional is
\[
\e^{-5} \left(\|1\|_{H^{-1}(T_\e\Gamma)}^2 - \frac23 \e^3\ell(\Gamma) - 2 n(\Gamma) \alpha \e^4 \right),
\]
which we again expect to approximate
\[
\frac2{45} \sum_i \ell(\gamma_i) \int_{\gamma_i} \kappa_i^2.
\]
The theory used in this paper to prove $\upGamma$-convergence is not adequately equipped to deal with open curves. For example, the notion of systems of curves includes only closed curves. An extension is needed to deal with the open curves.

Another, perhaps more approachable, question concerns the relation between $\|1\|_{H^{-1}(T_\e\gamma)}^2$ and $\|\chi_{T_\e\gamma}\|_{H^{-1}(\R^2)}^2$, where $\chi_{T_\e\gamma}$ is the characteristic function of the set $T_\e\gamma$. The latter expression is closer to what one can find in many applications, like the previously mentioned systems that exhibit partial localisation (Section~\ref{sec:motivation}).

Other extensions that bridge the gap between the current results and those applications a bit further are the study of $\|1\|_{H^{-1}(\Omega)}^2$ on neighbourhoods of curves that have a variable thickness or research into the $H^{-1}$-norm of more general functions, $\|f\|_{H^{-1}(T_\e\gamma)}^2$.

\subsection{Structure of the paper}

We start out in Section~\ref{sec:formalintroduction} with a formal calculation for closed curves which serves as a motivation for the results in Theorems~\ref{th:compactness} and~\ref{thm:gammaconvergence}. In Section~\ref{sec:systemsclosedcurves} we give the definitions of system of curves and various related concepts. In our computations we use a parametrisation of the $T_\e\gamma$ which is specified in Section~\ref{subsec:parametrisation}.
Section~\ref{sec:proofcomplb}  is then devoted to the proof of the compactness and $\upGamma$-convergence results (Theorems~\ref{th:compactness} and~\ref{thm:gammaconvergence}). In Section~\ref{sec:opencurves} we state and prove the asymptotic development~\pref{devel:main} for open curves (Theorem~\ref{th:devel}).

\section{A formal calculation}\label{sec:formalintroduction}

We now give some formal arguments to motivate the statements of our main results for closed curves, and also to illustrate some of the technical difficulties. In this description we restrict ourselves to a single, simple, smooth, closed curve $\gamma$. 

Since the definition of $\mathcal G_\e$ implies that the global radius of curvature $\rho(\gamma)$ is bounded from below by $\epsilon$, we can parametrise $T_\e\gamma$ in the obvious manner. We choose one coordinate, $s\in[0,1]$, along the curve and the other, $t\in[-1,1]$, in the direction of the normal to the curve. As we show in Lemma~\ref{lem:parametrisationclosedcurves}, this parametrisation leads to the following characterisation of the $H^{-1}$-norm:
\begin{align}
\notag
\|1\|_{H^{-1}(T_\e\gamma)}^2 &= \sup \Biggl\{\int_0^1 \int_{-1}^1 \Biggl( 2 f(s,t) \epsilon \ell(\gamma) \bigl(1-\epsilon t \kappa(s)\bigr) - \frac{\epsilon (f_{,s})^2(s,t)}{(1-\epsilon t \kappa(s)) \ell(\gamma)}+\\
&\hspace{2.8cm} - (f_{,t})^2(s,t) \left(\frac1{\epsilon} - t \kappa(s)\right) \ell(\gamma) \Biggr)\,dt\,ds\Biggr\},
\label{char:HMOf}
\end{align}
where the supremum is taken over functions $f\in W^{1,2}$ that satisfy $f(s, \pm 1)=0$, and subscripts~$,s$ and $,t$ denote differentiation with respect to $s$ and $t$.

The corresponding Euler-Lagrange equation is
\begin{equation}\label{eq:eulerlagrangeclosedcurves}
\epsilon \ell(\gamma) \bigl(1-\epsilon t \kappa(s)\bigr) + \epsilon \left( \frac{f_{,s}(s,t)}{\bigl(1-\epsilon t \kappa(s)\bigr) \ell(\gamma)}\right)_{\hspace{-0.2cm},s} + \Bigl(( f_{,t}(s,t) \bigl(\epsilon^{-1} - t\kappa(s) \bigr) \ell(\gamma) \Bigr)_{\hspace{-0.1cm},t} = 0.
\end{equation}
Formally we solve this equation by using an asymptotic expansion
\[
f(s,t) = f_0(s,t) + \epsilon f_1(s,t) + \epsilon^2 f_2(s,t) + \epsilon^3 f_3(s,t) + \epsilon^4 f_4(s,t) + \ldots
\]
as \emph{Ansatz}. The boundary condition $f(s, \pm 1)=0$ should be satisfied for each order of $\epsilon$ separately. Substituting this into (\ref{eq:eulerlagrangeclosedcurves}) and collecting terms of the same order in $\epsilon$ we find for the first five orders
\begin{align}
f_{0,tt}(s,t) &= 0 &\Longrightarrow f_0(s,t) &= 0,\nonumber\\
f_{1,tt}(s,t) &= 0 &\Longrightarrow f_1(s,t) &= 0,\nonumber\\
f_{2,tt}(s,t) &= -1 &\Longrightarrow f_2(s,t) &= \frac12 (1-t^2),\nonumber\\
f_{3,tt}(s,t) &= -t \kappa(s) &\Longrightarrow f_3(s,t) &= \frac16 t \kappa(s) (1-t^2),\nonumber\\
f_{4,tt}(s,t) &= -\frac16 \kappa^2(s) (9t^2-1) &\Longrightarrow f_4(s,t) &= \frac1{24} \kappa^2(s) (-3 t^4 + 2 t^2 + 1).\label{eq:formalAnsatz}
\end{align}
Note that this {Ansatz} is reasonable only for closed curves, since the ends of a tubular neighbourhood have different behaviour.
These orders suffice to compute the $H^{-1}$-norm up to order $\epsilon^5$:
\begin{equation}\label{eq:formalexpansion}
\|1\|_{H^{-1}(\Omega_{\epsilon})}^2 = \frac23 \ell(\gamma) \epsilon^3 + \frac2{45} \epsilon^5 \ell(\gamma) \int_\gamma \kappa^2 + \mathcal{O}(\epsilon^6) \quad \text{ as } \epsilon\to 0.
\end{equation}

For a \emph{fixed} curve $\gamma\in W^{2,2}$, this expansion can be made rigorous. Theorem~\ref{th:devel} proves an extended version~\pref{devel:main} of this development, in which ends are taken into account.

For a sequence of \emph{varying} curves $\gamma_n$, on the other hand, the explicit dependence of~$f_3$ on~$\kappa$ in this calculation is a complicating factor. Even if a sequence~$\gamma_n$ converges strongly in $W^{2,2}$---and that is a very strong requirement---then the associated curvatures $\kappa_n$ converge in $L^2$. There is no reason for the \emph{derivatives} $\kappa_n'(s)$ to remain bounded in $L^2$, and the same is true for the derivatives $f_{n3,s}$. Therefore the second term under the integral in~\pref{char:HMOf}, which is formally of order $O(\e^6)$, may turn out to be larger, and therefore interfere with the other orders. In Theorems~\ref{th:compactness} and~\ref{thm:gammaconvergence} this problem is addressed by introducing a regularized version of $\kappa$ in the definition of $f_3$.

\medskip

The formal calculation we did in this section suggests that we need information about the optimal function $f$ in (\ref{char:HMOf}) up to a level $\e^4$. However, as we will see, for the proof of the lower bound part of Theorem~\ref{thm:gammaconvergence} (part 1) it suffices to use information up to order $\e^3$, (\ref{eq:triallb}). The reason why becomes apparent if we look in more detail at the calculation that led to the formal expansion in (\ref{eq:formalexpansion}). The contributions to this expansion involving $f_4$ are given by
\[
2 \ell(\gamma) \int_0^1 \int_{-1}^1 \Bigl( f_4(s,t) - f_{2,t}(s,t) f_{4,t}(s,t)\Bigr)\,dt\,ds = \frac1{12} \ell(\gamma) \kappa^2 \int_0^1 \int_{-1}^1 \Bigl(-15 t^4 + 6 t^2 + 1\Bigr)\,dt\,ds = 0.
\]
This means that replacing $f_4$ by $\hat f_4\equiv0$ does not change the expansion up to order $\e^5$ given in (\ref{eq:formalexpansion}). It is an interesting question to ponder whether this is a peculiarity of the specific function under investigation or a symptom of a more generally valid property.

Note that for the proof of the upper bound statement of Theorem~\ref{thm:gammaconvergence} (part 2) we do need a trial function that has terms up to order $\e^4$, (\ref{eq:trialub}).

\section{Systems of closed curves}\label{sec:systemsclosedcurves}

From now on we aim for rigour. The first task is to carefully define systems of curves, their equivalence classes, and notions of convergence. We only consider closed curves, and systems of closed curves, and therefore we use the unit torus $\vz T = \vz R/\vz Z$ as the common domain of parametrisation.

Let $\vz{T}_{(i)}$ be disjoint copies of $\vz{T}$ and let
\[
\underset{i}{\coprod} \vz{T}_{(i)}:= \bigcup_i \bigl\{(s,i): s\in \vz{T}_{(i)} \bigr\}
\]
denote their disjoint union. A \emph{$W^{1,2}$-system of curves} is a map $\Gamma: \overset{m}{\underset{i=1}{\coprod}} \vz{T}_{(i)} \to \R^2$ given by
\[
\Gamma(s, i) = \gamma_i(s),
\]
where $m\in\vz{N}$ and, for all $1\leq i \leq m$, $\gamma_i \in W^{1,2}(\vz{T};\R^2)$ is a \emph{closed curve} \emph{parametrised proportional to arc length} (i.e. $|\gamma_i'|$ is constant). The \emph{number of curves} in $\Gamma$ is defined as $\#\Gamma := m$. 
We denote such a system by
\[
\Gamma = \{\gamma_i\}_{i=1}^m.
\]
Analogously we define a \emph{$W^{2,2}$-system of curves}.

A system $\Gamma$ is called \emph{disjoint} if for all $i\neq j$, $(\gamma_i)\cap(\gamma_j) = \emptyset$.
A $W^{2,2}$-system of curves is said to be \emph{without transverse crossings} if for all $i,j\in\{1,\ldots, m\}$ and
all $s_1, s_2\in\vz{T}$,
\begin{equation}
\label{cond:tangency}
\gamma_i(s_1)=\gamma_j(s_2) \Rightarrow \gamma_i'(s_1) = \pm \gamma_j'(s_2).
\end{equation}
The \emph{length} of a curve $\gamma$ and of a system  of curves $\Gamma$ is
\[
\ell(\gamma) := \int_\T |\gamma'|
\qquad\text{and}\qquad \ell(\Gamma) := \sum_{i=1}^m \ell(\gamma_i).
\]
The global radius of curvature of a system of curves $\Gamma$ is 
\[
\rho(\Gamma) := \underset{x, y, z \in (\Gamma)}{\inf} r(x, y, z),
\]
where $r(x, y, z)$ is the radius of the unique circle in $\R^2$ through $x$, $y$, and $z$ if $x, y, z\in (\Gamma)$ are pairwise disjoint and not collinear and $r(x,y,z)=\infty$ otherwise, analogous to Definition~\ref{def:globalcurvature}. The $\e$-tubular neighbourhood of $\Gamma$ is the set $T_\e\Gamma$,
\[
T_\e\Gamma:= \bigcup_{x\in(\Gamma)} B(x,\e),
\]
where $B(x, \e)$ denotes the open ball with center $x$ and radius $\e$.

Let $\{\Gamma^n\}_{n=1}^{\infty}$ be a sequence of $W^{k,2}$-systems of curves, $k=1, 2$. We write $\Gamma^n = \{\gamma_i^n\}_{i=1}^m$ and say $\Gamma^n$ \emph{converges} to $\Gamma$ \emph{in $W^{k,2}$} for a $W^{k,2}$-system of curves $\Gamma=\{\gamma_i\}_{i=1}^m$ if for $n$ large enough $\#\Gamma^n = \#\Gamma$ and for all $1\leq i \leq m$, $\gamma^n_i \to \gamma_i$ in $W^{k,2}(\vz{T};\R^2)$ as $n\to\infty$ (after reordering). We write $\Gamma^n \to \Gamma$ in $W^{k,2}$. 

The \emph{density function} $\theta_{\Gamma}: (\Gamma) \to \vz{N}\cup\{+\infty\}$ of a system of curves $\Gamma$ is defined as
\[
\theta_{\Gamma}(z) := \mathcal{H}^0(\{\Gamma^{-1}(z)\}).
\]

Let $\Gamma$ and $\tilde \Gamma$ be two $W^{2,2}$-systems of curves.
We say that $\Gamma$ and $\tilde \Gamma$ are \emph{equivalent}, denoted by $\Gamma \sim \tilde \Gamma$, if $(\Gamma)=(\tilde\Gamma)$ and $\theta_{\Gamma}=\theta_{\tilde\Gamma}$ everywhere. We denote the set of \emph{equivalence classes of $W^{k,2}$-systems of curves},
$k\in\{1,2\}$, by $SC^{k,2}$. Where necessary we explicitly write $[\Gamma]$ for the equivalence class that contains $\Gamma$; where possible we will simply write $\Gamma$ for both the system of curves and for its equivalence class.

Let $\{[\Gamma^n]\}_{n=1}^{\infty},\ [\Gamma]\subset SC^{k,2}$, $k\in\{1,2\}$. We say that $[\Gamma^n]$ \emph{converges} to $[\Gamma]$
in $SC^{k,2}$  if there exist $\Gamma^n\in[\Gamma^n]$ and $\Gamma\in[\Gamma]$ such that $\Gamma^n \to \Gamma$ in $W^{k,2}$ in the sense defined above. We denote this convergence by $[\Gamma^n]\to [\Gamma]$ in $SC^{k,2}$.

\begin{rmrk}
\label{rem:indepofrep}
Note that if $\tilde \Gamma\in [\Gamma]$, then $(\Gamma)=(\tilde\Gamma)$, so that the definition $([\Gamma]):=(\Gamma)$ is independent of the choice of representative.  Similarly, the length $\ell(\Gamma)$, the curvature $\kappa$, the global radius of curvature $\rho(\Gamma)$, the tubular neighbourhood $T_\e\Gamma$, the functional $\mathcal G_\e$, and the property of having no transverse crossings are all well-defined on equivalence classes. The same is also true for the functional $\mathcal G_0$; this is proved in~\cite[Lemma 3.9]{BellettiniMugnai04}.
\end{rmrk}

\begin{rmrk}
If $\gamma\in W^{k,2}(\vz{T};\R^2)$, $k=1,2$, is a curve parametrised proportional to arc length, then it follows that
\[
|\gamma'| = \ell(\gamma) \qquad \text{(a.e. if }k=1).
\]
\end{rmrk}

We also introduce some elementary geometric notation. 
Let $\gamma \in W^{2,2}\left(\vz{T}; \R^2\right)$ be
a curve parametrised proportional to arc length.
We choose the {normal} to the curve at $\gamma(s)$ to be
\[
\nu(s) := |R \gamma'(s)|^{-1}R \gamma'(s)=\ell(\gamma)^{-1} R \gamma'(s),
\]
where $R$ is the anticlockwise rotation matrix given by
\[
R := \left( \begin{array}{ll} 0&-1\\1&0 \end{array} \right).
\]
The \emph{curvature} $\kappa: \T \to \R$ satisfies
\begin{equation}
\label{def:curvature}
\gamma''(s) = \kappa(s) \ell(\gamma)^2 \nu(s).
\end{equation}
We have
\begin{align*}
|\nu|=1, &\qquad \nu' = \kappa \ell(\gamma) R \nu,\\
\gamma' \times \nu = \ell(\gamma), &\qquad \nu' \times \nu = -\kappa \ell(\gamma),
\end{align*}
where $\times$ denotes the cross product in $\R^2$:
\[
x \times y := x_1 y_2 - x_2 y_1 = (Rx)\cdot y, \qquad \text{ for } x, y \in \R^2.
\]
It is well known that integrating the curvature of a closed curve gives
\begin{equation}\label{eq:integratingthecurvature}
\ell(\gamma) \int_\T \kappa = -\int_\T \nu'\times\nu = \pm 2 \pi,
\end{equation}
depending on the direction of parametrisation. Without loss of generality we adopt a parametrisation convention which gives the $+$-sign in the integration above, and which could be described as `counterclockwise'. The integral of the squared curvature can be expressed as
\[
\int_\gamma \kappa^2 = \ell(\gamma)\int_\T \kappa(s)^2\, ds 
= \ell(\gamma)^{-3} \int_\T |\gamma''(s)|^2 \, ds.
\]

\section{Parametrising the tubular neighbourhood}\label{sec:setup}
\label{subsec:parametrisation}

By density we have 
\[
\|1\|_{H^{-1}(T_\e\gamma)}^2 = \sup\left\{ \int_{T_\e\gamma} \left( 2 \phi(x) - |\nabla \phi(x)|^2\right)\,dx : \phi\in W_0^{1,2}(T_\e\gamma) \right\}, 
\]
and the supremum is achieved when $\phi$ equals $\varphi\in C^\infty(\T_\e\gamma)\cap C(\overline{T_\e\gamma})$, the solution of
\begin{equation}\label{eq:thisiseulerlagrangeforZ}
\left\{ \begin{array}{ll}
  -\Delta\varphi = 1&\text{ in } T_\e\gamma,\\ 
  \varphi = 0 & \text{ on } \partial T_\e\gamma. 
\end{array}\right.
\end{equation}
In that case we also have
\[\|1\|_{H^{-1}(T_\e\gamma)}^2 = \int_{T_\e\gamma} |\nabla \varphi(x)|^2\,dx
\]

In the proof of our main result, Theorem~\ref{thm:gammaconvergence}, we use a reparametrisation of the $\epsilon$-tubular neighbourhood of a simple $W^{2,2}$-closed curve. For easy reference we introduce it here in a separate lemma.

\begin{figure}[h]
\centering
    \psfrag{b}{\small{$\gamma(s)-\epsilon \nu(s)$}}
    \psfrag{c}{\small{$\gamma'(s)$}}
    \psfrag{d}{\small{$\gamma(s)+\epsilon \nu(s)$}}
    \psfrag{e}{\small{$\nu'(s)$}}
    \includegraphics[width=50mm]{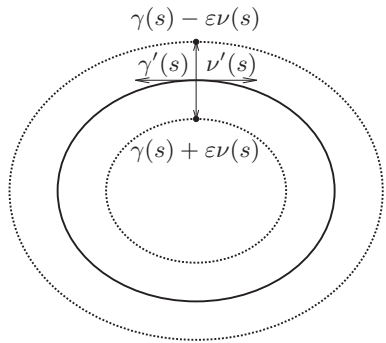}\\
\caption{A closed curve with an $\epsilon$-tubular neighbourhood. Explicitly shown are normal $\nu(s)$ and tangent $\gamma'(s)$ at $\gamma(s)$ and the points $\gamma(s)\pm\epsilon \nu(s)$}\label{fig:closedcurve}
\end{figure}

\begin{lmm}\label{lem:parametrisationclosedcurves}
Let $\epsilon>0$ and let $\gamma \in W^{2,2}(\vz{T};\R^2)$ be a
closed curve parametrised proportional to arc length, and such that $\rho(\gamma)\geq \epsilon$. If we define $\Psi_\e \in W^{1,2}\left(\vz{T} \times (-1,1);T_\e\gamma\right)$ by
\begin{equation}\label{eq:parametrisation}
\Psi_\e(s,t) := \gamma(s) + \epsilon t \nu(s),
\end{equation}
then $\Psi_\e$ is a bijection.

Let $g \in W^{1,2}\left(T_\e\gamma\right)$ and define $f := g \circ \Psi_\e$. Then
\[
\int_{T_\e\gamma} \left(2 g(x) - |\nabla g(x)|^2\right)\,dx = \mathcal{X}_\e(f),
\]
where
\begin{align*}
\mathcal{X}_\e(f) &:= \int_0^1 \int_{-1}^1 \Biggl( 2 f(s,t) \epsilon \ell(\gamma) \bigl(1-\epsilon t \kappa(s)\bigr) - \frac{\epsilon (f_{,s})^2(s,t)}{\bigl(1-\epsilon t \kappa(s)\bigr) \ell(\gamma)}+\\&\hspace{5cm} - (f_{,t})^2(s,t) \left(\frac1{\epsilon} - t \kappa(s)\right) \ell(\gamma) \Biggr)\,dt\,ds.
\end{align*}
\end{lmm}
The parametrisation of $T_\e\gamma$ from Lemma~\ref{lem:parametrisationclosedcurves} is illustrated in Figure~\ref{fig:closedcurve}.

\begin{proof}[Proof of Lemma~\ref{lem:parametrisationclosedcurves}]
We first show that
$\Psi_\e: \vz{T} \times (-1, 1) \to T_\e\gamma$ is a bijection. Starting with surjectivity, we fix $x\in T_\e\gamma$; by the discussion in Section~\ref{sec:intro-curvature} there exists a unique $s\in\vz{T}$ such that $\gamma(s)$ is the point of minimal distance to $x$ among all points in $(\gamma)$. The line segment connecting $x$ and $\gamma(s)$ necessarily intersects $(\gamma)$ perpendicularly and thus there exists a $t \in (-1,1)$ such that $x=\Psi_\e(s,t)$.

We prove injectivity by contradiction. Assume there exist
$(s, t), (\tilde s, \tilde t)\in \vz{T} \times (-1, 1)$ and $x\in T_\e\gamma$, such that $(s, t) \neq (\tilde s, \tilde t)$ and $\Psi_\e(s,t)=\Psi_\e(\tilde s, \tilde t) = x$. If $s=\tilde s$, then $t \neq \tilde t$, which contradicts $\Psi_\e(s,t)=\Psi_\e(\tilde s, \tilde t)$, so we assume now that $s\neq \tilde s$. Also without loss of generality we take $\tilde t \leq t < 1$.
We compute
\begin{equation}\label{eq:gammagamma}
\gamma(\tilde s)-\gamma(s) = \e \bigl(t \nu(s) - \tilde t \nu(\tilde s)\bigr).
\end{equation}
Let $r(\gamma(\tilde s), \gamma(s), z)$ be as in Definition~\ref{def:globalcurvature} and let $\theta$ be the angle between $\gamma(\tilde s)-\gamma(s)$ and $\gamma'(\tilde s)$. By \cite[Equation 3]{GonzalezMaddocks99} if we take the limit $z\to \gamma(s)$ along the curve we find
\begin{align*}
r(\gamma(\tilde s), \gamma(s), \gamma(s)) &= \frac{|\gamma(s)-\gamma(\tilde s)|}{2 |\sin\theta|}
= \frac{\ell(\gamma) |\gamma(\tilde s)-\gamma(s)|^2}{2 |\gamma'(s) \times (\gamma(\tilde s)-\gamma(s))|}\\
&= \frac{\e^2 \ell(\gamma) |\tilde t \nu(\tilde s) - t \nu(s)|^2}{2 \e |\gamma'(s) \times (\tilde t \nu(\tilde s) - t \nu(s))|}
= \e \frac{t^2 + \tilde t^2 - 2 t \tilde t \nu(s) \cdot \nu(\tilde s)}{2 |\tilde t \nu(s)\cdot \nu(\tilde s)-t|}.
\end{align*}
Note that $\tilde t \nu(s)\cdot\nu(\tilde s) \leq \tilde t \leq t$ and
\[
t^2+\tilde t^2 - 2 t \tilde t \nu(s)\cdot\nu(\tilde s) - 2\bigl(t-\tilde t \nu(s)\cdot\nu(\tilde s)\bigr) = 2\bigl(t-\tilde t \nu(s)\cdot\nu(\tilde s)\bigr) (t-1) + \tilde t^2 - t^2 < 0,
\]
from which we conclude that $r(\gamma(\tilde s), \gamma(s), \gamma(s)) < \e$ which contradicts $\rho(\gamma)\geq \e$. Therefore, $\Psi_\e$ is injective and thus a bijection.

We compute
\[
\nabla f^T = \left(\nabla g \circ \Psi_\e\right)^T D\Psi_\e,
\]
where $D\Psi_\e$ is the derivative matrix of $\Psi_\e$ in the $(s,t)$-coordinates. It follows that
\[
|\nabla g|^2 \circ \Psi_\e = \nabla f^T D\Psi_\e^{-1} D\Psi_\e^{-T} \nabla f,
\]
where
$\cdot^{-T}$ denotes the inverse of the transpose of a matrix. Direct computation yields
\[
D\Psi_\e(s,t) = \left( \begin{array}{ll} 
\gamma_1'(s) - \epsilon \ell(\gamma) t \kappa(s) \nu_2(s)
     & \epsilon \nu_1(s)\\
\gamma_2'(s) + \epsilon \ell(\gamma) t \kappa(s) \nu_1(s)
     &\epsilon \nu_2(s) 
\end{array} \right)
\]
and $\det D\Psi_\e(s,t) = \epsilon \ell(\gamma) \left(1-\epsilon t \kappa(s)\right)$. Since
$\|\kappa\|_{L^{\infty}(\vz{T})}\leq \epsilon^{-1}$  we have $\det D\Psi_\e(s,t) \neq 0$ almost everywhere. Then
\[
(D\Psi_\e)^{-1}(s,t) (D\Psi_\e)^{-T}(s,t) = \left(\begin{array}{ll} \ell(\gamma)^{-2} \left(1-\epsilon t \kappa(s)\right)^{-2}&0\\0&\epsilon^{-2}\end{array}\right), 
\]
and we compute
\begin{align*}
&\hspace{0.5cm}\int_{T_\e\gamma} \left(2 g(x) - |\nabla g(x)|^2\right)\,dx\\
&= \int_0^1 \int_{-1}^1 \left( 2 f(s,t) - \frac{(f_{,s})^2(s,t)}{\ell(\gamma)^2 \bigl(1-\epsilon t \kappa(s)\bigr)^2} - \frac{(f_{,t})^2(s,t)}{\epsilon^2} \right) |\det D\Psi_\e(s,t)|\,dt\,ds,
\end{align*}
which gives the desired result.
\end{proof}

The previous lemma gives us all the information to compute the $H^{-1}$-norm of $1$ on a tubular neighbourhood:

\begin{crllr}\label{cor:H-1innewparametrisation}
Let $\gamma\in W^{2,2}(\vz{T};\R^2)$ be a
closed curve parametrised proportional to arc length with $\rho(\gamma) \geq \epsilon$. Furthermore let $\Psi_\e, \mathcal{X}_\e$ be as in Lemma~\ref{lem:parametrisationclosedcurves}. Define
\begin{equation}\label{eq:setofadmissibles}
\mathcal{A}_{\epsilon} := \left\{f\in W^{1,2}\left(\vz{T}\times[-1,1]\right) : f\circ\Psi_\e^{-1}\in W_0^{1,2}\left(T_\e\gamma\right)\right\}.
\end{equation}
Then
\begin{equation}\label{eq:gradensupzaken}
\|1\|_{H^{-1}(T_\e\gamma)}^2 = \sup\left\{\mathcal{X}_\e(f): f\in\mathcal{A}_{\epsilon}\right\}.
\end{equation}
\end{crllr}

\drop{
The next lemma shows that $\mathcal{G}_{\epsilon}$ is nonnegative.

\begin{lmm}
Let $\epsilon>0$. For all $[\Gamma]\in SC^{1,2}$, $\mathcal{G}_{\epsilon}([\Gamma]) \geq 0$.
\end{lmm}
\begin{proof}
If $[\Gamma]\notin\mathcal{S}_{\epsilon}$ the statement follows trivially, so assume that $[\Gamma]\in\mathcal{S}_{\epsilon}$. Let $\gamma\in\Gamma\in[\Gamma]$, let $\Omega_{\epsilon}$ be the $\epsilon$-tubular neighbourhood around $(\gamma)$, and define $f\in W_0^{1,2}(\vz{T}\times[-1,1])$ by $f(s,t):= \frac12 \epsilon^2 (1-t^2)$.
Then in the notation of Lemma~\ref{lem:parametrisationclosedcurves} and using Corollary~\ref{cor:rewriteH-1norm} we find that
\[
\|1\|_{H^{-1}(\Omega_{\epsilon})} \geq \mathcal{X}(f) = \frac23 \epsilon^3 \ell(\gamma).
\]
Now use Remark~\ref{rem:indepofrepres} and in particular (\ref{eq:sumalot1}) to sum over all curves.
\end{proof}
}

\section{Proof of Theorem~\ref{th:compactness} and the lower bound part of Theorem~\ref{thm:gammaconvergence}}\label{sec:proofcomplb}

\subsection{Reduction to single curves}

Let us first make a general remark. If $\mathcal G_\e(\Gamma)$ is finite, then $\rho(\Gamma)\geq \e$, and therefore the $\epsilon$-tubular neighbourhoods of two distinct curves in $\Gamma$ do not intersect. Therefore writing $\Gamma=\{\gamma_n\}_{n=1}^m$, we can decompose $\mathcal G_\e(\Gamma)$ as
\begin{equation}\label{eq:sumalot1}
\mathcal{G}_{\epsilon}(\Gamma) = \sum_{n=1}^m \mathcal{G}_{\epsilon}(\gamma_n).
\end{equation}
A similar property also holds for $\mathcal{G}_0$ if $\Gamma\in\mathcal{S}_0$, as follows directly from the definition:
\begin{equation}\label{eq:sumalot2}
\mathcal{G}_0(\Gamma) = \sum_{n=1}^m \mathcal{G}_0(\gamma_n).
\end{equation}

\subsection{Trial function}
\label{sec:trialfunction}

The central tool in the proof of compactness (Theorem~\ref{th:compactness}) and the lower bound inequality (part \ref{thm:gammaconvergence:1} of Theorem~\ref{thm:gammaconvergence}) is the use of a specific choice of $f$ in $\mathcal{X}_\e(f)$. For a given $\gamma\in W^{2,2}(\T,\R^2)$, this trial function is of the form
\begin{equation}\label{eq:triallb}
f_\e(s,t) = \frac{\e^2}2(1-t^2) + {\e^3}\bar \kappa_\e(s)\zeta(t).
\end{equation}
Here $\bar \kappa_\e$ is an $\e$-dependent approximation of $\kappa$ which we specify in a moment, and $\zeta\in C^1_c(-1,1)$ is a fixed, nonzero, odd function satisfying
\begin{equation}
\label{cond:zeta}
\int_{-1}^1 {\zeta'}^2(t)\,dt = \int_{-1}^1 t\zeta(t)\,dt .
\end{equation}
In the final stage of the proof $\zeta$ will be chosen to be an approximation of the function $t(1-t^2)/6$. Note that this choice for $f$ can be seen as an approximation of the first two non-zero terms in the asymptotic development~\pref{eq:formalAnsatz}.  As explained at the end of Section~\ref{sec:formalintroduction} this suffices and we do not need a term of order $\e^4$ in $f$.

When used in $\mathcal X_\e$, the even and odd symmetry properties in $t$ of the two terms in $f_\e$ cause various terms to cancel. The result is
\[
\|1\|^2_{H^{-1}(T_\e\gamma)}\geq \mathcal{X}_\e(f_{\epsilon}) = \frac23 \epsilon^3 \ell(\gamma) 
+ B\e^5\ell(\gamma)\int_\T \Bigl\{2  \kappa(s) \bar\kappa_{\epsilon}(s) - \bar\kappa_{\epsilon}^2(s) - \e^2\tilde C_{\epsilon}(s) \bar\kappa_{\epsilon}'^2(s)\Bigr\}\,ds,
\]
where
\begin{align}
B &:= \int_{-1}^1 \zeta'^2(t)\,dt = \int_{-1}^1 t\zeta(t)\,dt ,\label{eq:Bepsn}\\
\tilde C_{\epsilon}(s) &:= B^{-1}\ell(\gamma)^{-2} \int_{-1}^1 \frac{\zeta^2(t)}{(1-\epsilon t \kappa(s))}\,dt,\nonumber
\end{align}
The definition of $\tilde C_\e$ shows why $\zeta$ is chosen with compact support in $(-1,1)$. By the uniform bound $\|\kappa\|_\infty \leq \epsilon^{-1}$, the denominator $1-\e t\kappa(s)$ is uniformly bounded away from zero, \emph{independently of the curve $\gamma$}. Therefore $\ell(\gamma)^2\tilde C_\e$ is bounded from above and away from zero
independently of $\gamma$.

It will be convenient to replace the $s$-dependent coefficient $\tilde C_\e$ by a constant coefficient. For that reason we introduce
\[
C_\e := \sup_{s\in\T} \tilde C_\e(s),
\]
which is finite for fixed $\e$ and $\gamma$.
With this we have 
\[
\|1\|^2_{H^{-1}(T_\e\gamma)}\geq  \frac23 \epsilon^3 \ell(\gamma) 
+ B\e^5\ell(\gamma)\int_\T \Bigl\{2 \kappa(s) \bar\kappa_{\epsilon}(s) -\bar\kappa_{\epsilon}^2(s) - \e^2C_{\epsilon} \bar\kappa_{\epsilon}'^2(s)\Bigr\}\,ds.
\]
This expression suggests a specific choice for $\bar \kappa_\e$: choose $\bar\kappa_\e$ such as to maximize the expression on the right-hand side. The Euler-Lagrange equation for this maximization reads
\begin{equation}
\label{eq:EulerKappa}
- \e^2C_\e\bar\kappa_\e''(s) + \bar\kappa_\e(s) = \kappa(s)
\qquad\text{for a.e. }s\in \T,
\end{equation}
from which the regularity $\bar\kappa_\e\in W^{2,2}(\T)$ can be directly deduced; this regularity is sufficient to guarantee $f_\e\circ \Psi_\e^{-1}\in W_0^{1,2}(T_\e\gamma)$, so that the resulting function $f_\e$ is admissible in $\mathcal A_\e$ (see~\pref{eq:setofadmissibles}). 
The resulting maximal value provides the inequality
\begin{equation}\label{eq:estimatefixedepsilongamma}
\|1\|^2_{H^{-1}(T_\e\gamma)}\geq  \frac23 \epsilon^3 \ell(\gamma) 
+ B\e^5\ell(\gamma)\int_\T \Bigl\{ \bar\kappa_{\epsilon}^2(s) + \e^2C_{\epsilon} \bar\kappa_{\epsilon}'^2(s)\Bigr\}\,ds.
\end{equation}

\subsection{Step 1: fixed number of curves}

We now place ourselves in the context of Theorem~\ref{th:compactness}. Let $\epsilon_n\to 0$  and $\{\Gamma^n\}_{n=1}^{\infty}\subset SC^{1,2}$ be sequences such that
$\mathcal{G}_{\epsilon_n}(\Gamma^n)$ and $\ell(\Gamma^n)$ are bounded uniformly by a constant $\Lambda>0$. We need to prove that there exists a subsequence of the sequence $\{\Gamma^n\}$ that converges in $SC^{1,2}$ to a limit $\Gamma\in \mathcal{S}_0$. 

The first step is to limit the analysis to a fixed number of curves, 
which is justified by the following lemma.

\begin{lmm}\label{lem:compactness-boundednumber}
There exists a constant $C>0$ depending only on $\Lambda$ such that 
\[
\frac1C\leq \ell(\gamma)\leq  C \qquad\text{for any } n\in \N \text{ and any }\gamma\in \Gamma_n.
\]
Consequently $\#\Gamma_n$ is bounded uniformly in $n$.
\end{lmm}
\begin{proof}[Proof of Lemma~\ref{lem:compactness-boundednumber}]
For any $n$ choose an arbitrary $\gamma\in\Gamma^n$; then $\mathcal{G}_{\epsilon_n}(\gamma)\leq \Lambda$, and therefore
by~\pref{eq:estimatefixedepsilongamma}  the associated $\bar\kappa_{\e_n}$ satisfies
\[
\int_\T \bar \kappa_{\e_n}^2 \leq \frac \Lambda{B\ell(\gamma)}.
\]
Integrating~\pref{eq:EulerKappa} over $\T$ and using periodicity we then find
\begin{equation}
\label{eq:intbarkappaepsilonis2piab-1}
2\pi = \ell(\gamma) \int_\T \kappa
  = \ell(\gamma) \int_\T \bar\kappa_{\e_n}
  \leq \ell(\gamma) \Bigl(\int_\T \bar\kappa_{\e_n}^2 \Bigr)^{1/2}
  \leq \ell(\gamma) \left(\frac \Lambda{B\ell(\gamma)}\right)^{1/2},
\end{equation}
which implies that $\ell(\gamma)$ is bounded from below; therefore any curve in any $\Gamma^n$ has its length bounded from below. Since $\ell(\Gamma^n)$ is bounded from above, the result holds.
\end{proof}

Because of this result, we can restrict ourselves to a subsequence along which $\#\Gamma^n$ is constant. We switch to this subsequence without changing notation. 

Moreover, by the discussion in Section~\ref{sec:trialfunction} we find that $C_{\e_n}$ is bounded uniformly in $\e_n$ and $\gamma$. We can therefore apply inequality~\pref{eq:estimatefixedepsilongamma} to any sequence of curves $\gamma_n$, corresponding to a sequence $\e_n\to0$, as in the statements of Theorems~\ref{th:compactness} and~\ref{thm:gammaconvergence}. In the terminology of those theorems we find the inequality
\begin{equation}
\liminf_{n\to\infty} \mathcal{G}_{\epsilon_n}(\gamma_n) 
\geq   \liminf_{n\to\infty} B\ell(\gamma_n) \int_\T \Bigl\{ \bar\kappa_{\epsilon_n}^2(s) + \e_n^2C_{\epsilon_n} \bar\kappa_{\epsilon_n}'^2(s)\Bigr\}\,ds.
\label{eq:liminfbounded2}
\end{equation}

\subsection{Step 2: single-curve analysis}

For every $n$, we pick an arbitrary curve $\gamma\in\Gamma^n$, and for the rest of this section we label this curve $\gamma_n$. The aim of this section is to prove appropriate compactness properties and the lower bound inequality for this sequence of single curves. 

In this section we associate with the sequence $\{\gamma_n\}$ of curves the curvatures $\kappa_n$ (see~\pref{def:curvature}) and the quantities $\bar\kappa_n:=\bar\kappa_{\e_n}$ and $C_n:= C_{\e_n}$ that were introduced in Section~\ref{sec:trialfunction}.
Note that by~\pref{eq:liminfbounded2}, the upper bound on $\mathcal{G}_{\e_n}(\gamma_n)$, and the lower bound on $\ell(\gamma_n)$ there exists an $M>0$ such that 
\begin{equation}
\label{eq:energyboundbarkappa2}
\int_\T \Bigl\{ \bar\kappa_n^2(s) + \e_n^2C_n \bar\kappa_n'^2(s)\Bigr\}\,ds \leq M 
\qquad\text{for all }n\in\N.
\end{equation}

\begin{lmm}\label{lem:compactness-limitsofcurvatures}
There exists a subsequence of $\{\bar\kappa_n\}_{n=1}^{\infty}$ (which we again label by $n$), such that
\begin{equation}\label{eq:weakL2convergenceofbarkappa}
\bar\kappa_n\longrightharpoonup \bar \kappa \text{ in } L^2(\vz{T}),
\end{equation}
for some $\bar\kappa \in L^2(\vz{T})$, and
\begin{equation}\label{eq:kappaepsionndistributionalconvergence}
\kappa_n \longrightharpoonup \bar\kappa \text{ in } H^{-1}(\vz{T}).
\end{equation}
In addition, defining
\begin{equation}
\label{def:theta_n0}
\vartheta_n(s) := \ell(\gamma_n) \int_0^s \kappa_n(\sigma)\,d\sigma
\qquad\text{and}\qquad
\vartheta_0(s) := \ell(\gamma_0) \int_0^s \bar\kappa(\sigma)\,d\sigma,
\end{equation}
we have
\[
\vartheta_n \weakto \vartheta_0 \qquad\text{in }L^2(\T).
\]
\end{lmm}
\begin{proof}[Proof of Lemma~\ref{lem:compactness-limitsofcurvatures}]
By (\ref{eq:energyboundbarkappa2}), $\{\bar\kappa_n\}_{n=1}^{\infty}$ is uniformly bounded in $L^{2}(\vz{T})$, and therefore there is a subsequence (which we again index by $n$) such that $\bar \kappa_n \longrightharpoonup \bar \kappa$ in $L^2(\vz{T})$ for some $\bar \kappa$
in $L^2(\vz{T})$.

Next let $f\in C^1(\vz{T})$ and compute
\begin{align}
\int_\T \kappa_n(s)f(s)\,ds &= \int_\T \left(\bar\kappa_n(s) - C_n \epsilon_n^2  \bar\kappa_n''(s)\right)f(s)\,ds\nonumber\\
&= \int_\T \left(\bar\kappa_n(s)f(s) + C_n\epsilon_n^2 \bar\kappa_n'(s) f'(s)\right)\,ds.
\label{eq:distributionalderivative}
\end{align}
By the uniform lower bound on $\ell(\gamma_n)$ (Lemma~\ref{lem:compactness-boundednumber}) and (\ref{eq:energyboundbarkappa2}) we have for some $C\geq 0$
\[
\int_\T |\bar\kappa_n'(s) f'(s)| \,ds \leq \|\bar\kappa_n'\|_{L^2(\vz{T})} \|f'\|_{L^2(\vz{T})} \leq C \epsilon_n^{-1}.
\]
Therefore the last term in (\ref{eq:distributionalderivative}) converges to zero and thus $\kappa_n$ converges weakly to
$\bar\kappa$ in~$H^{-1}(\vz{T})$. From the definition~\pref{def:theta_n0} and the convergence $\ell(\gamma_n)\to\ell(\gamma_0)$ it then follows that $\vartheta_n\weakto \vartheta_0$ in $L^2(\T)$.
\end{proof}

We next bootstrap the weak $L^2$-convergence of $\vartheta_n$ to strong $L^2$-convergence.

\begin{lmm}
After extracting another subsequence (again without changing notation) we have
\[
\vartheta_n\to\vartheta_0\qquad\text{in }L^2(\T)\text{ and pointwise a.e.}
\]
\end{lmm}

\begin{proof}
For the length of this proof it is more convenient to think of all functions as defined on $[0,1]$ rather than on $\T$. Define $\bar K_n\in W^{1,2}(0,1)$
by
\[
\bar K_n(s) := \ell(\gamma_n) \int_0^s \bar\kappa_n(t)\,dt.
\]
By the bound on $\ell(\gamma_n)$ in Lemma~\ref{lem:compactness-boundednumber} we can set for the duration of this proof $\ell(\gamma_n)=1$ without loss of generality.
The boundedness of $\bar\kappa_n=\bar K_n'$ in $L^2(0,1)$ (see~\pref{eq:energyboundbarkappa2}) implies that $\bar K_n$ is compact in $C^{0,\alpha}([0,1])$ for all $0< \alpha<1/2$. By integrating (\ref{eq:EulerKappa}) from $0$ to
$s>0$
we find
\[
\vartheta_n(s) - C_n \epsilon_n^2  \bar \kappa_n'(0) = \bar K_n(s) - C_n \epsilon_n^2  \bar \kappa_n'(s).
\]
Inequality (\ref{eq:energyboundbarkappa2}) also gives
\[
C_n\int_0^1 \epsilon_n^4 \bar\kappa_n'^2(s)\,ds \leq \epsilon_n^2 M,
\]
by which
\[
C_n \epsilon_n^2\bar \kappa_n' \to 0 \quad \text{ in } L^2(0,1),
\]
and combined with the compactness of $\bar K_n$ in $C^{0,\alpha}([0,1])$ this implies that 
$\left\{\vartheta_n(s) - C_n \epsilon_n^2 \bar \kappa_n'(0)\right\}_{n=1}^{\infty}$ is compact in
$L^2(0,1)$. Since we already know that $\vartheta_n$ converges weakly in $L^2(0,1)$, it follows that (along a subsequence) the sequence of constant functions $C_n \epsilon_n^2  \bar \kappa_n'(0)$ converges weakly, i.e. that the scalar sequence $C_n \epsilon_n^2  \bar \kappa_n'(0)$ converges in $\R$. Therefore $\vartheta_n$ converges strongly to~$\vartheta_0$.

\end{proof}

Let us write
\[
\gamma'_n(s) = \ell(\gamma_n) \left(\begin{array}{c} \cos(\vartheta_n(s)+\varphi_n)\\\sin(\vartheta_n(s)+\varphi_n)\end{array}\right),
\]
where $\varphi_n\in [0,2\pi)$ is an $n$-dependent phase.

We then use the uniform boundedness of $\gamma_n(0)\in B(0,R)$ and of $\varphi_n\in[0,2\pi)$ to extract yet another subsequence such that $\gamma_n(0)$ converges to some $x_0\in \overline {B(0,R)}$ and $\varphi_n$ converges to some $\varphi\in[0,2\pi]$. Defining the curve $\gamma_0$ by
\begin{equation}
\label{param:theta_0}
\gamma_0(0) := x_0 
\qquad\text{and}\qquad
\gamma'_0(s) := \ell(\gamma_0) \left(\begin{array}{c} \cos(\vartheta_0(s)+\varphi)\\\sin(\vartheta_0(s)+\varphi)\end{array}\right),
\end{equation}
it follows from the strong convergence of $\vartheta_n$ in $L^2(\T)$ that $\gamma_n\to\gamma_0$ in the strong topology of $W^{1,2}(\T;\R^2)$.

\medskip
We can now find an $L^2$-bound on $\gamma_0''$.

\begin{lmm}\label{lem:compactness-boundoncurvature}
We have 
\[
\|\gamma_0''\|_{L^2(\vz{T};\R^2)} = \ell(\gamma_0)^2 \|\bar \kappa\|_{L^2(\vz{T})}.
\]
\end{lmm}
\begin{proof}
Since $\vartheta_0'= \ell(\gamma_0)\bar\kappa\in L^2(\T)$, upon differentiating~\pref{param:theta_0} we find
\[
\gamma_0''(s) = \ell(\gamma_0) \vartheta_0'(s)
  \left(\begin{array}{c} -\sin(\vartheta_0(s)+\varphi)\\\cos(\vartheta_0(s)+\varphi)\end{array}\right),
  \qquad\text{at a.e. }s\in \T,
\]
so that 
\[
\|\gamma_0''\|_{L^2(\T;\R^2)} = \ell(\gamma_0) \|\vartheta_0'\|_{L^2(\T)}
= \ell(\gamma_0)^2\|\bar\kappa\|_{L^2(\T)}.
\]

\end{proof}

\subsection{Step 3: Returning to systems of curves}

We have shown that the sequence of single curves $\{\gamma_n\}$ satisfies $\gamma_n\to\gamma_0$ in $W^{1,2}(\T;\R^2)$ with $\ell(\gamma_0)<\infty$ and
$\|\gamma_0''\|_{L^2(\vz{T}; \R^2)}= \ell(\gamma_0)^2\|\bar\kappa\|_{L^2(\T)}<\infty$. For future reference we note that this implies that
\begin{equation}
\label{ineq:liminf_prep}
\mathcal G_0(\gamma_0) = \frac 2{45}\ell(\gamma_0)^{-3} \|\gamma_0''\|_{L^2(\T;\R^2)}^2
= \frac 2{45}\ell(\gamma_0)\|\bar\kappa\|_{L^2(\T)}^2
\leq \frac2{45B} \liminf_{n\to\infty} \mathcal G_{\e_n} (\gamma_n).
\end{equation}
The inequality follows from (\ref{eq:liminfbounded2}), (\ref{eq:weakL2convergenceofbarkappa}), and the weak-lower semicontinuity of the $L^2$-norm.

Now we return from the sequence of single curves to the sequence of systems of curves $\{\Gamma^n\}_{n=1}^{\infty}$. Write $\Gamma^n := \{\gamma_i^n\}_{i=1}^m$, and repeat the above arguments for each sequence of curves $\{\gamma_i^n\}_{n=1}^{\infty}$ for fixed $i$ separately. In this way we find a limit system $\Gamma^0:=\{\gamma_i^0\}_{i=1}^m$ such that for all $i$, $\ell(\gamma_i^0)<\infty$ and $\|(\gamma_i^0)''\|_{L^2(\vz{T}; \R^2)}<\infty$.
It is left to prove that $\Gamma^0$ has no transverse crossings.

\begin{lmm}\label{lem:compactness-notransverse crossings}
$\Gamma_0$ has no transverse crossings.
\end{lmm}
\begin{proof}[Proof of Lemma~\ref{lem:compactness-notransverse crossings}]
We prove this by contradiction.

Assume that $\Gamma^0$ has transverse crossings. This can happen if either two different curves in $\Gamma^0$ intersect transversally or if one curve self-intersects transversally. First assume the former, i.e. assume that there exist $\gamma_1^0, \gamma_2^0 \in \Gamma^0$ and $s_1, s_2\in\vz{T}$ such that $\gamma_1^0(s_1)=\gamma_2^0(s_2)$ and
${\gamma_1^0}'(s_1) \neq \pm {\gamma_2^0}'(s_2)$. Without loss of generality we take $s_1=s_2=0 \in \vz{T}$ and $\gamma_1^0(0)=0$. For ease of notation in this proof we will identify $\vz{T}$ with the interval $[-1/2,1/2]$ with the endpoints identified.

Because $\gamma_1^0, \gamma_2^0 \in C^1(\vz{T};\R^2)$ and ${\gamma_1^0}'(s_1) \neq \pm {\gamma_2^0}'(s_2)$ there exists a $\delta>0$ such that 
\[
\text{if }s,t\in[-\delta,\delta]\text{ satisfy }\gamma_1^0(s)=\gamma_2^0(t), \text{ then }s = t=0. 
\]

Define
\[
D := (-\delta, \delta) \times (-\delta, \delta) \subset \R^2
\]
and the function $f\in C^1(\overline{D}; \R^2)$
 by
\[
f(s,t) := \gamma_1^0(s) - \gamma_2^0(t).
\]
We compute
\[
\det Df(s,t) = {\gamma_2^0}'(t) \times {\gamma_1^0}'(t)
\]
and find that
\[
\det Df(0,0) \neq 0,
\]
since we assumed that ${\gamma_1^0}'$ and ${\gamma_2^0}'$ are not parallel. Since $f(s,t)=0$ iff $(s,t)=(0,0)$ and furthermore $0\not\in f(\partial D)$ we can use \cite[Definition 1.2]{FonsecaGangbo95} to compute the {topological degree} of $f$  with respect to $D$:
\[
d(f, D, 0) := \sum_{(s,t)\in f^{-1}(0)} \text{sgn} \left(\det Df(s,t)\right)  =  \text{sgn} \left(\det Df(0,0)\right) = \pm 1,
\]
where the sign depends on the direction of parametrisation of $\gamma_1^0$ and $\gamma_2^0$.

We know that $\Gamma^n \to \Gamma^0$ in $W^{1,2}$ as $n\to\infty$, so in particular for $n$ large enough there are curves $\gamma_1^n, \gamma_2^n \in \Gamma^n$ such that
\[
\gamma_i^n \to \gamma_i^0 \quad \text{ in } C(\vz{T}; \R^2) \text{ as } n\to\infty, \,i\in\{1,2\}.
\]

If we now define $f^n\in C^1(\overline{D}; \R^2)$ by
\[
f^n(s,t) := \gamma_1^n(s) - \gamma_2^n(t),
\]
then we conclude by \cite[Theorem 2.3 (1)]{FonsecaGangbo95} that for large enough $n$,
\[
d(f^n, D, 0) = d(f, D, 0) \neq 0.
\]

By  \cite[Theorem 2.1]{FonsecaGangbo95},  $d(f^n, D, 0) \neq 0$ implies that there exists $(s_0,t_0)\in D$ such that $f^n(s_0,t_0) = 0$, i.e. that $\gamma_1^n(s_0) = \gamma_2^n(t_0)$. This contradicts the fact that $(\gamma_1^n) \cap (\gamma_2^n) = \emptyset$ and therefore we deduce that $\Gamma^0$ does not contain two different curves that cross transversally.

Now assume that a single curve $\gamma^0\in \Gamma^0$ has a transversal self-intersection. Then we can repeat the above argument with $\gamma_1^0=\gamma_2^0=\gamma^0$ and $s_1\neq s_2$ to deduce that there exist $\gamma^n\in \Gamma^n$ such that $\gamma_n^0 \to \gamma^0$ in $C(\T; \R^2)$ and  $\gamma^n(s_0) = \gamma^n(t_0)$, which contradicts the bound on the global curvature of $\Gamma^n$.

We conclude that $\Gamma^0$ has no transverse crossings.
\end{proof}

This concludes the proof of the compactness result from Theorem~\ref{thm:gammaconvergence}.

\subsection{Proof of the lower bound of Theorem~\ref{thm:gammaconvergence}}

Let $[\Gamma^n]\to[\Gamma]$ be a sequence as in part~\ref{thm:gammaconvergence:1} of Theorem~\ref{thm:gammaconvergence}. Then by the  definition of the convergence of equivalence classes, we can choose representatives $\tilde\Gamma^n\in[\Gamma^n]$ and $\tilde\Gamma\in[\Gamma]$ such that $\tilde\Gamma^n\to\tilde \Gamma$. We drop the tildes for notational convenience. By the definition of convergence of systems of curves, for $n$ large enough $\Gamma^n=\{\gamma_i^n\}_{i=1}^m$, i.e. the number $m$ of curves is fixed, and each $\Gamma^n$ can be reordered such that $\gamma_i^n \to \gamma_i$ in $W^{1,2}(\vz{T}; \R^2)$ for each $i$.

Without loss of generality we assume that $\liminf_{n\to\infty} \mathcal G_{\e_n}(\Gamma^n)< \infty$, and that for all $N$,
\[
\mathcal G_{\e_N^{}}(\Gamma^N)\leq \liminf_{n\to\infty} \mathcal G_{\e_n}(\Gamma^n) +1.
\]
Since $W^{1,2}(\T;\R^2)\subset L^\infty(\T;\R^2)$ and $\gamma_i^n \to \gamma_i$ in $W^{1,2}(\vz{T}; \R^2)$, the traces $(\gamma_i^n)$ are all contained in some large bounded set. 
Therefore Theorem~\ref{th:compactness} applies, and there exists a subsequence along which $\gamma_i^n$ converges in $W^{1,2}(\T;\R^2)$ to limit curves $\gamma_i^0$. Since limits are unique, we have $\gamma_i^0 = \gamma_i$.

We then calculate by~\pref{ineq:liminf_prep}
\[
\mathcal G_0(\Gamma) = \sum_{i=1}^m  \mathcal G_0(\gamma_i)
\leq \frac2{45B} \sum_{i=1}^m \liminf_{n\to\infty} \mathcal G_{\e_n}(\gamma_i^n)
= \frac2{45B} \liminf_{n\to\infty} \mathcal G_{\e_n}(\Gamma^n).
\]
The required liminf bound follows by remarking that by choosing $\zeta\in C_c^\infty(-1,1)$ odd, satisfying~\pref{cond:zeta}, and close to the function $\tilde \zeta(t) := t(1-t^2)/6$, the number $B$ can be chosen arbitrarily close to $2/45$.

This concludes the proof of part~\ref{thm:gammaconvergence:1} of Theorem~\ref{thm:gammaconvergence}.

\subsection{Proof of the $\limsup$ inequality from Theorem~\ref{thm:gammaconvergence}}\label{sec:limsupproof}

For a single, fixed, simple, smooth, closed curve $\gamma$, the formal calculation of Section~\ref{sec:formalintroduction} can be made rigorous. This is done in the context of open curves in Lemma~\ref{lemma:bulk}, and the argument there can immediately be transferred to closed curves. For such a curve therefore
\[
\lim_{n\to\infty}\mathcal G_{\e_{n}}(\gamma) = \mathcal G_{0}(\gamma).
\]

The only remaining issue is therefore to show that any~$\Gamma$ can be approximated by a system~$\tilde\Gamma$ consisting of smooth, disjoint, simple closed curves. This is the content of the following lemma. 

\begin{lmm}\label{lem:detach}
Let $\Gamma$ be a $W^{2,2}$-system of closed curves without
transversal crossings. Then there exists a number $m>0$, a sequence of systems $ \left\{\Gamma^j\right\}_{j=1}^\infty$, and a system $\tilde \Gamma = \{\tilde \gamma_k\}_{k=1}^m$ equivalent to $\Gamma$ such that the following holds:
\begin{enumerate}
\item
For all $j\in\N$ the system of curves
$
  \Gamma^j = \{\gamma^j_k\}_{k=1}^m
$
is a pairwise disjoint family of smooth simple closed curves;
\item
For all $1\leq k\leq m$ we have
\begin{align}\label{eq:conv-strong-gamma}
  & \gamma^j_k \,\to\, \tilde \gamma_k\quad\text{ in }W^{2,2}(0,1)\quad\text{ as }
    j\to\infty.
\end{align}
\end{enumerate}
In particular we have $[\Gamma^j]\to [\Gamma]$ in $SC^{2,2}$ and 
$  \mathcal G_0([\Gamma^j])\to \mathcal G_0([\Gamma])$ as $j\to\infty$.
\end{lmm}

This lemma is proved in~\cite[Lemma~8.2]{PeletierRoeger08}. By taking a diagonal sequence the lim sup inequality follows.

\section{Open curves}\label{sec:opencurves}

The aim of this section is to prove the following theorem:
\begin{thrm}
\label{th:devel}
Let $\gamma\in C^\infty([0,1];\R^2)$ satisfy
\begin{itemize}
\item $\gamma$ is parametrised proportional to arclength; 
\item $\gamma$ is exactly straight (i.e.\ $\gamma''\equiv0$) on a neighbourhood of each end.
\end{itemize}
Then there exists a constant $\alpha>0$ (see~\pref{def:alpha}), independent of $\gamma$, such that
\begin{equation}
\label{devel:main_v2}
\|1\|_{H^{-1}(T_\e\gamma)}^2 = \frac23 \e^3 \ell(\gamma) + 2\alpha \e^4 + 
  \frac 2{45}\e^5\int_{\gamma} \kappa^2 + O(\e^6)\qquad\text{as }\e\to0.
\end{equation}
\end{thrm}

\subsection{Overview of the proof}
\begin{figure}[h]
    \centering\let\tiny\footnotesize
    \psfrag{a}{\tiny$\e$}
    \psfrag{b}{$\Omega_{\eta}$}
    \psfrag{c}{$\Omega_{\eta}$}
    \psfrag{d}{$\Omega_-$}
    \psfrag{e}{\tiny$\upGamma$}
    \psfrag{f}{\tiny$\upGamma$}
    \psfrag{g}{\tiny$\gamma(1)$}
    \psfrag{h}{\tiny$\gamma(0)$}
    \psfrag{i}{\tiny$2\eta$}
    \psfrag{j}{\tiny$2\eta$}
    \includegraphics[width=0.35\textwidth]{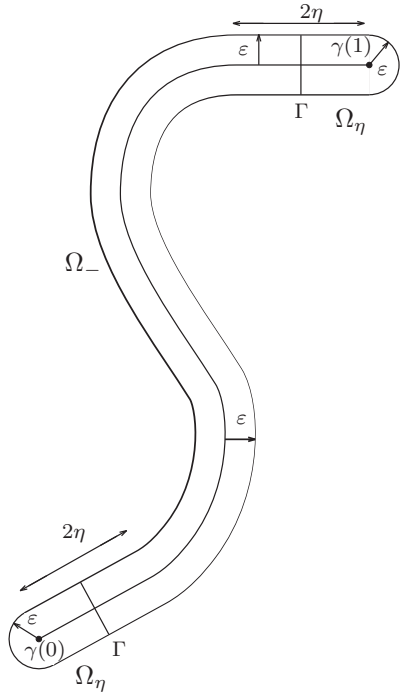}\\
    \caption{$\e$-tubular neighbourhood of an open curve with two straight endings}\label{fig:opencurvetwosidesflat}
\end{figure}

The proof of Theorem~\ref{th:devel} hinges on a division of the domain into separate parts. 
To make this precise we introduce some notation. 

First we note that the squared $H^{-1}$-norm in two dimensions scales as $(\mathrm{length})^4$, i.e. if $\Omega\subset \R^2$, then 
\[
\|1\|_{H^{-1}(\lambda\Omega)}^2 = \lambda^4 \|1\|_{H^{-1}(\Omega)}^2.
\]
Therefore the development~\pref{devel:main_v2} is scale-invariant under a rescaling of both $\gamma$ and~$\e$ by a common factor (i.e. a rescaling of $T_\e\gamma$ by this same factor); by multiplying both by $\ell(\gamma)^{-1}$ we can assume that the curve $\gamma$ has length $1$.

Next we define the normal $\nu$ and the curvature $\kappa$ as in Section~\ref{sec:systemsclosedcurves}. We also use the parametrisation
\[
\Psi_\e(s,t) := \gamma(s) + \e t\nu(s), 
\]
although for an open curve $\Psi_\e$ only covers the tubular neighbourhood without the end caps. 

We let $0<2\eta<1$ be a length of parametrisation corresponding to the straight end sections, i.e. we choose $\eta$ such that
\[
\gamma''(s) = 0 \qquad\text{for }s\in [0,2\eta]\cup [1-2\eta,1].
\]
We then define 
\[
\Omega_- := \Psi_\e\bigl((\eta,1-\eta)\times (-1,1)\bigr).
\]
Note that $\Omega_-$ contains the bulk of the tubular neighbourhood, and \emph{half} of each of the straight sections near the ends. The remainder, corresponding to two end caps with the other half of the straight sections, is
\[
\Omega_\eta := T_\e\gamma \setminus \overline{\Omega_-}.
\]
We call $\upGamma = \overline {\Omega_-}\cap \overline{\Omega_\eta}$ the interface separating $\Omega_-$ from $\Omega_\eta$. 
See Figure~\ref{fig:opencurvetwosidesflat}.
\medskip
The statement of Theorem~\ref{th:devel} follows from the following three lemmas. The first implies that we may cut up the domain $T_\e\gamma$ into $\Omega_-$ and $\Omega_\eta$ and consider the two domains separately.

\begin{lmm}
\label{lemma:devel:cutting}
Define the boundary data function $\ubc :\upGamma \to\R$ by
\[
\ubc(\Psi_\e(s,t)) = \frac{\e^2}2 (1-t^2)\quad\text{for }s\in \{\eta,1-\eta\}, \ t\in(-1,1).
\]
Then 
\[
\|1\|_{H^{-1}(T_\e\gamma)}^2 = 
\int_{\Omega_-} u_- + \int_{\Omega_\eta} u_\eta + O(e^{-\eta/\e}), \quad \text{as } \e\to 0,
\]
where $u_-:\Omega_-\to\R$ and $u_\eta : \Omega_\eta\to\R$ are solutions of
\begin{equation}
\label{def:u-ueta}
\left\{\begin{aligned}
-\Delta u_- &= 1&\quad&\text{in }\Omega_-\\
u_- &= 0 &\quad &\text{on }\partial\Omega_-\setminus\upGamma\\
u_- &= \ubc &\quad&\text{on }\upGamma
\end{aligned}\right.
\qquad\qquad
\left\{
\begin{aligned}
-\Delta u_\eta &= 1&\quad&\text{in }\Omega_\eta\\
u_\eta &= 0 &\quad &\text{on }\partial\Omega_\eta\setminus\upGamma\\
u_\eta &= \ubc &\quad&\text{on }\upGamma.
\end{aligned}\right.
\end{equation}
\end{lmm}

The second lemma deals with the bulk of the tubular neighbourhood.

\begin{lmm}
\label{lemma:bulk}
We have
\[
\int_{\Omega_-} u_- = \frac23\e^3(1-2\eta)
 + \frac 2{45}\e^5\int_{\gamma} \kappa^2 + O(\e^6), \quad \text{as } \e\to 0.
\]
\end{lmm}

The third lemma gives an estimate of the contribution of the ends. 
\begin{lmm}
\label{lemma:ends}
We have
\[
\int_{\Omega_\eta} u_\eta = \frac43\eta \e^3 + 2\alpha \e^4 + O(e^{-\eta/\e}), \quad \text{as } \e\to 0,
\]
where $\alpha>0$ is given in~\pref{def:alpha}.
\end{lmm}

\subsection{Proof of Lemma~\ref{lemma:devel:cutting}}
By partial integration we have
\[
\|1\|_{H^{-1}(T_\e\gamma)}^2 = \int_{T_\e\gamma} u,
\]
where $u: T_\e\gamma\to\R$ solves
\[
-\Delta u = 1\quad \text{in }T_\e\gamma, 
\qquad
u = 0 \quad\text{on }\partial T_\e\gamma.
\]

We first note the useful property that there exists a constant $M$, independent of $\e$, such that 
\[
\|u\|_{L^\infty(T_\e\gamma)}\leq M.
\]
This follows from remarking that all $T_\e\gamma$ are contained in a large ball $B(0,R)$, and that the solution of 
\[
-\Delta v = 1 \quad \text{in }B(0,R), 
\qquad
v = 0 \quad\text{on }\partial B(0,R)
\]
is a supersolution for $u$, independent of $\e$. Without loss of generality we can assume that $M\geq 1$. 

We next turn to the content of the lemma. Below we show that
\begin{equation}
\label{est:lemma:cutting:Linfty}
\|u_--u\|_{L^\infty(\Omega_-)} + \|u_\eta-u\|_{L^\infty(\Omega_\eta)} = O(e^{-\eta/\e}),
\end{equation}
from which the assertion follows, since
\begin{align*}
\left|\int_{T_\e\gamma} u - \int_{\Omega_-} u_- - \int_{\Omega_\eta} u_\eta\right|
&\leq \int_{\Omega_-}|u-u_-| + \int_{\Omega_\eta} |u-u_\eta|\\
&\leq |\Omega_-|\|u-u_-\|_{L^\infty(\Omega_-)}  + |\Omega_\eta|\|u_\eta-u\|_{L^\infty(\Omega_\eta)}.
\end{align*}

To show~\pref{est:lemma:cutting:Linfty} we first consider an auxiliary problem, that we formulate as a lemma for future reference. 
\begin{lmm}
\label{lem:maxprinc2}
Define the rectangle $\mathcal R$ and its boundary parts, 
\begin{align*}
\mathcal{R} &:= (-a,a) \times (-b,b),\\
\partial \mathcal{R}_1 &:= \left\{ (x,y)\in\partial\mathcal{R} : |y|=b \right\},\\
\partial \mathcal{R}_2 &:= \left\{ (x,y)\in\partial\mathcal{R} : |x|=a \right\}
\end{align*}
and let $g\in C^{\infty}(\mathcal{R})\cap C(\overline{\mathcal R})$ satisfy
\[
\left\{\begin{array}{ll} -\Delta g = 0 & \text{ on } \mathcal{R},\\ g=0 & \text{ on } \partial\mathcal{R}_1,\\ |g|\leq 1 & \text{ on } \partial\mathcal{R}_2. \end{array}\right.
\]
Then 
\[
|g(0,y)| \leq 4e^{-a/b} \quad \text{for all }y\in(-b,b).
\]
\end{lmm}
The proof of this lemma follows from remarking that 
\[
\chi(x,y) := \frac{\cosh({x}/b) \cos({y}/b)}{\cosh(a/b)\cos 1}
\]
is a supersolution for this problem, and therefore 
\[
|g(0,y)| \leq \chi(0,y) \leq 4e^{-a/b} \quad \text{for all }y\in(-b,b).
\]

We now apply this estimate to the straight sections at each end of $\gamma$. Assume that one of the straight sections coincides with the rectangle $\mathcal R$ with $a=\eta$ and $b=\e$ (this amounts to a translation and rotation of $\gamma$). Note that then the line segment $\{0\}\times(-\e,\e)$ is part of $\upGamma$.  The function $g(x,y) := \tfrac12M^{-1}[u(x,y) - \frac12 (\e^2-y^2)]$ satisfies the conditions above, and therefore
\[
|u(0,y) - \tfrac12 (\e^2-y^2)| = | u(0,y) - \ubc(0,y)| =  O(e^{-\eta/\e}), \qquad\text{uniformly in }y\in (-\e,\e).
\]
At the other end a similar estimate holds, implying that
\[
\|u-\ubc\|_{L^\infty(\upGamma)} =  O(e^{-\eta/\e}).
\]
We then deduce the estimate~\pref{est:lemma:cutting:Linfty} by applying the maximum principle to $u-u_-$ in $\Omega_-$ and to $u-u_\eta$ in $\Omega_\eta$.

\subsection{Proof of Lemma~\ref{lemma:bulk}}

As in Section~\ref{subsec:parametrisation} we can write
\[
\int_{\Omega_-} u_- = \e\int_\eta^{1-\eta} \int_{-1}^1
  u_-(\Psi_\e(s,t))\,(1-\e t\kappa(s))\, dsdt.
\]
For the length of this section we set $\omega := (\eta,1-\eta) \times (-1,1)$. 
Writing $f_-(s,t) := u_-(\Psi_\e(s,t))$, we find
\begin{equation}\label{eq:twicethesame}
\int_{\Omega_-} \bigl( 2u_- - |\nabla u_-|^2 \bigr) = \int_\omega \bigl( 2 f_- \e (1-\e t \kappa) - \nabla f_- \cdot B_\e \nabla f_-\bigr),
\end{equation}
where
\[
B_{\epsilon}(t,s) := \left( 
  \begin{array}{ll} 
    \epsilon  \bigl(1-\epsilon t \kappa(s)\bigr)^{-1} & 0 \\ 
    0 & \epsilon^{-1} \bigl(1-\epsilon t \kappa(s)\bigr) 
  \end{array}
\right).
\]
By \pref{def:u-ueta} $u_-$ satisfies the Euler-Lagrange equation corresponding to the left hand side of (\ref{eq:twicethesame}). Therefore $f_-$ satisfies the Euler-Lagrange equation for the right hand side:
\[
\left\{ \begin{array}{ll}
-\div  B_{\epsilon} \nabla f_- = \e (1-\e t \kappa) & \text{on } \omega,\\ 
f_-(s, \pm 1) = 0 & \text{for }s\in (\eta,1-\eta),\\ 
f_-(s,t) = \frac{\e^2}2(1-t^2) & \text{for }(s,t)\in \{\eta,1-\eta\}\times (-1,1).
\end{array}\right.
\]

We now define the trial function
\begin{equation}\label{eq:trialub}
f_\e (s,t) := \frac{\e^2}2(1-t^2) + \frac{\e^3}6\kappa(s) t(1-t^2) 
  + \frac{\e^4}{24}\kappa^2(s)(-3t^4+2t^2+1)
\end{equation}
for which we calculate that
\begin{alignat*}2
\div  B_{\epsilon} \nabla (f_--f_\e) &= h_\e &\qquad& \text{in } \omega,\\ 
(f_--f_\e)(s, \pm 1) &= 0, && \text{on }\partial \omega,
\end{alignat*}
where the defect $h_\e$ satisfies
\[
\|h_\e\|_{L^\infty(\omega)} = O(\e^4).
\]

Below we prove that this estimate on $h_\e$ implies that
\begin{equation}
\label{est:f-fe}
\|f_--f_\e\|_{L^2(\omega)} = O(\e^5).
\end{equation}
Assuming this estimate for the moment, we
find the statement of the lemma by the same calculation as in Section~\ref{sec:formalintroduction},
\[
\e\int_\eta^{1-\eta} \int_{-1}^1
  f_\e(s,t)\,(1-\e t\kappa(s))\, dsdt
  = \frac23\e^3(1-2\eta)
 + \frac 2{45}\e^5\int_{\gamma} \kappa^2 + O(\e^6),
\]
and the remark that
\[
\left|\e\int_\eta^{1-\eta} \int_{-1}^1
  (f_--f_\e)(s,t)\,(1-\e t\kappa(s))\, dsdt\right|
\leq 4(1-2\eta)\e\|f_--f_\e\|_{L^2((\eta,1-\eta)\times(-1,1))}.
\]

To prove~\pref{est:f-fe} we set $g=f_--f_\e$ and apply Poincar\'e's inequality in the $t$-direction:
\[
\int_{-1}^1 g^2(s,t)\,dt \leq \frac4{\pi^2} \int_{-1}^1 (g_{,t})^2(s,t)\,dt,
\]
where $g_{,t}$ again denotes the partial derivative of $g$ with respect to $t$. 
Since
\[
\epsilon^{-1}  \left(1-\epsilon t \kappa(s)\right) \left(g_{,t}(s,t)\right)^2 \leq \nabla g(s,t) \cdot B_{\epsilon}(s,t) \nabla g(s,t),
\]
we then calculate
\begin{align*}
\int_\eta^{1-\eta}\int_{-1}^1g^2(s,t)\, dtds 
&\leq \frac4{\pi^2} \int_\eta^{1-\eta}\int_{-1}^1 (g_{,t})^2(s,t)\,dtds\\
&\leq \frac4{\pi^2}(1+O(\e))\int_\eta^{1-\eta}\int_{-1}^1 (g_{,t})^2(s,t) \, (1-\e t\kappa(s))\,dtds\\
&\leq \frac{4\e}{\pi^2}(1+O(\e))\int_\eta^{1-\eta}\int_{-1}^1 
  \nabla g(s,t) \cdot B_{\epsilon}(s,t) \nabla g(s,t)\,dtds\\
&= -\frac{4\e}{\pi^2}(1+O(\e))\int_\eta^{1-\eta}\int_{-1}^1 
  g(s,t) \div B_{\epsilon}(s,t) \nabla g(s,t)\,dtds\\
&\leq \frac{4\e}{\pi^2}(1+O(\e)) \|g\|_{L^2(\omega)}\|h_\e\|_{L^2(\Omega)},
\end{align*}
so that 
\[
\|g\|_{L^2(\omega)} \leq \frac{4\e}{\pi^2}(1+O(\e)) \|h_\e\|_{L^2(\Omega)}
= O(\e^5).
\]

\subsection{Proof of Lemma~\ref{lemma:ends}}
\label{subsec:ends}

The domain $\Omega_\eta$ consists of two unconnected parts. We prove the result for just one of them, the part at the end $\gamma(1)$. We assume without loss of generality that 
\[
\gamma(1) = 0
\qquad\text{and}\qquad
\gamma\left([1-2 \eta, 1]\right) = \left\{(x,0)\in\R^2: -2 \eta \leq x \leq 0\right\}.
\]
Then the 
corresponding end of $\Omega_\eta$ is a reduction by a factor $\e$ of the domain
\[
\omega_{\eta}:= (-\eta/\e,0)\times(-1,1) \cup B(0,1)
\]
which itself is a truncation of the set
\[
\omega := (-\infty,0)\times(-1,1) \cup B(0,1).
\]
The sets $\omega_\eta$ and $\omega$ are depicted in Figure~\ref{fig:twoomegas}.

\def\scale{0.8}
\begin{figure}[htb]
\subfloat[Omegainfty][The domain $\omega$ in Section~\ref{subsec:ends}]
{
    \psfrag{a}[][][\scale]{$\omega$}
    \psfrag{b}[][][\scale]{$1$}
    \psfrag{c}[][][\scale]{$1$}
    \psfrag{d}[][][\scale]{$(0,0)$}
    \psfrag{e}[][][\scale]{$\leftarrow -\infty$}
    \includegraphics[height=1.8cm]{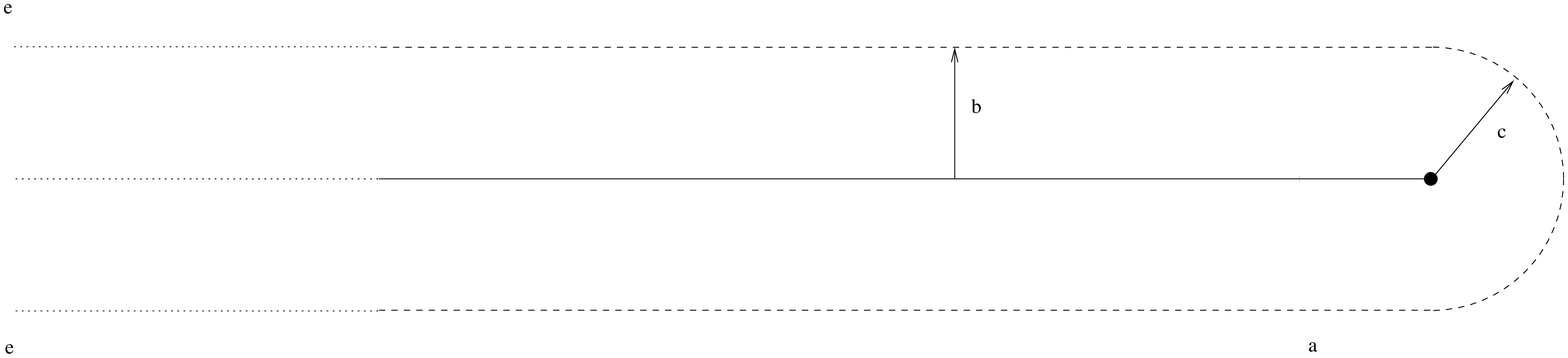}
    \label{fig:omegainfty}
}
\hspace{3mm}
\subfloat[Omegaeta][The domain $\omega_\eta$ in Section~\ref{subsec:ends}]
{
    \psfrag{a}[][][\scale]{$\omega_\eta$}
    \psfrag{b}[][][\scale]{$1$}
    \psfrag{c}[][][\scale]{$1$}
    \psfrag{d}[][][\scale]{$(0,0)$}
    \psfrag{e}[l][][\scale]{$(-\eta/\e, 0)$}
    \includegraphics[height=1.6cm]{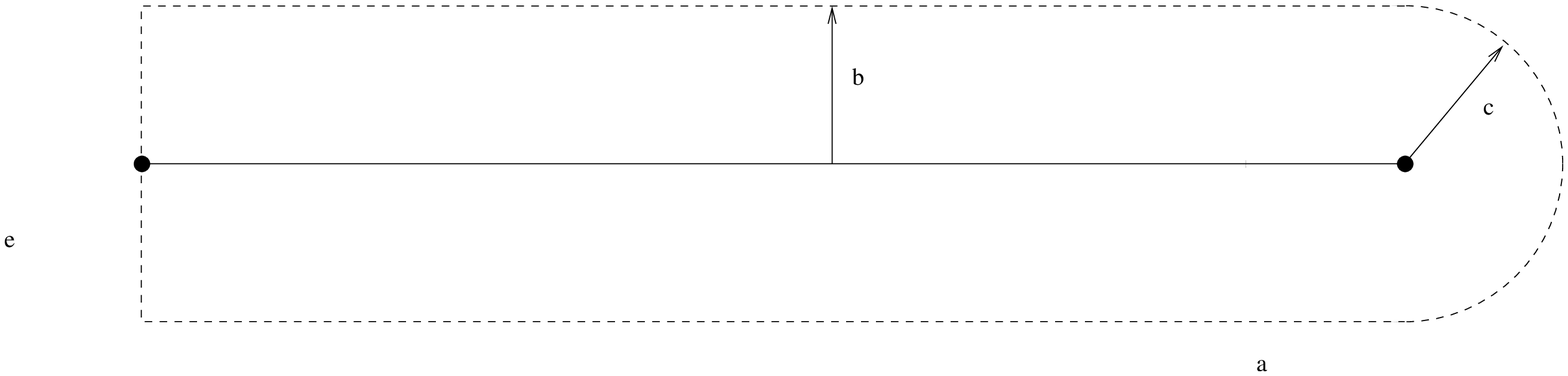}
    \label{fig:omegaeta}
}
\caption{}\label{fig:twoomegas}
\end{figure}

We also set $v_{\eta}(x,y) := \e^{-2}u_{\eta}(\e x,\e y)$, so that
\[
\int_{\Omega_{\eta}} u_{\eta} = \e^{4}\int_{\omega_{\eta}} v_{\eta}.
\]

Set 
\[
\varphi(x,y) := \frac12(1-y^2).
\]
The function $\psi_{\eta}:= v_{\eta}-\varphi$ then satisfies
\[
\left\{
\begin{aligned}
&{-\Delta\psi_{\eta}} = 0&\qquad&\text{in }\omega_{\eta},\\
&\psi_{\eta}=-\varphi&&\text{on }\partial\omega_{\eta}\setminus\{-\eta/\e\}\times(-1,1),\\
&\psi_{\eta}=0&&\text{on }\partial\omega_{\eta}\cap\{-\eta/\e\}\times(-1,1),
\end{aligned}
\right.
\]
and we have
\[
\e^{-4}\int_{\Omega_{\eta}} u_{\eta} = \int _{\omega_{\eta}} v_{\eta}
= \int_{\omega_{\eta}}\varphi + \int_{\omega_{\eta}} \psi_{\eta}.
\]
The first integral on the right-hand side is easily calculated:
\[
\int_{\omega_{\eta}}\varphi = \frac23\,\frac\eta\e + \frac3{16}\pi.
\]

Since we are taking the limit $\e\to0$, in which the set $\omega_{\eta}$ converges to the set $\omega$, we also define
$\psi\in C^{\infty}\left(\omega\right) \cap C\left(\overline\omega\right)$ to be the unique solution of
\[
\left\{
\begin{aligned}
&{-\Delta\psi} = 0&\qquad&\text{in }\omega,\\
&\psi=-\varphi&&\text{on }\partial\omega,\\\
&\|\psi\|_{L^{\infty}(\omega)} <\infty,
\end{aligned}
\right.
\]
and the constant 
\begin{equation}
\label{def:alpha}
\alpha := \int_{\omega} \psi + \frac{3}{16}\pi.
\end{equation}
Note that by the maximum principle $\|\psi\|_{L^{\infty}(\omega)} \leq \|\psi\|_{L^{\infty}(\partial\omega)}=1/2$.

Applying Lemma~\ref{lem:maxprinc2} to the rectangle $(2x,0)\times(-1,1)\subset\omega$ (with $x<0$) we find the decay estimate
\begin{equation}
\label{est:decay_psi}
|\psi(x,y)|\leq 2 e^{-x}\qquad\text{for all }y\in (-1,1)\text{ and all }x<0.
\end{equation}
This implies that
\[
\int _{\omega\setminus\omega_{\eta}}|\psi |\leq 4e^{-\eta/\e}.
\]
This estimate also provides an estimate of $\psi-\psi_{\eta}$. Note that $\psi=\psi_{\eta}=0$ on all of $\partial \omega_{\eta}$ with the exception of $\{-\eta/\e\}\times(-1,1)$. Applying~\pref{est:decay_psi} to this latter set we find that
\[
|\psi-\psi_\eta|  \leq 2 e^{-\eta/\e}\qquad \text{on all of }\partial \omega_{\eta}
\]
and since $\psi-\psi_{\eta}$ is harmonic on $\omega_{\eta}$ we conclude by the maximum principle that 
\[
\|\psi-\psi_{\eta}\|_{L^{\infty}(\omega_{\eta})} \leq 2e^{-\eta/\e}.
\]

The statement of Lemma~\ref{lemma:ends} then follows by remarking that
\[
\int_{\omega_{\eta}} u_{\eta} = \frac23 \eta\e^3 + \frac3{16}\pi\e^{4} 
+ \e^{4}\int _{\omega}\psi
+ R,
\]
where the rest term $R$ satisfies
\[
\e^{-4}|R| = \left|\int_{\omega_{\eta}} \psi_{\eta}-\int _{\omega}\psi\right|
\leq \int_{\omega_{\eta}}|\psi_{\eta}-\psi| + \int _{\omega\setminus\omega_{\eta}}|\psi|
\leq  2\left(\frac{2\eta}\e + \frac{\pi}2\right) e^{-\eta/\e} + 2e^{-\eta/\e}.
\]

\medskip
The only remaining assertion of the lemma is that $\alpha>0$. A finite-element calculation provides the estimate
\[
\alpha \approx 0.139917;
\]
here we only prove that $\alpha>0$. Define the harmonic comparison function
\begin{align*}
\tilde \psi(x,y) &= -0.112\cdot e^{\pi x/2}\cos (\pi y/2)
 +0.0019\cdot e^{3\pi x/2}\cos (3\pi y/2)
 -0.00008\cdot e^{5\pi x/2}\cos (5\pi y/2)\\
 &\hspace{0.45cm}- 0.056\cdot e^x\cos y,
\end{align*}
for which we can calculate (partially by numerical approximation of the appropriate one-dimensional integral)
\[
\int_{\omega}\tilde \psi \approx -0.5875 >  \frac3{16}\pi.
\]
We have $\psi\geq \tilde \psi$ on $\partial \omega$, implying that
\[
\alpha \;=\; \frac3{16}\pi + \int _{\omega} \psi \;\geq \;\frac3{16}\pi + \int _{\omega}\tilde \psi 
\;>\;0.
\]

%%-----------------------------
%%      your bibliography
%%-----------------------------

\bibliographystyle{siam}
%\bibliography{bibliography}
\bibliography{H-1_COCV-August14_2009.bbl}

\end{document}